\documentclass{amsproc}

\usepackage{pb-diagram}
\usepackage{amssymb}

\def\A{\mathcal {A}}

\def\S{\mathcal {S}}

\def\Fe{\mathcal {F}}
\def\D{\mathcal {D}}
\def\R{\mathbf{R}}
\def\C{\mathbf{C}}
\def\F{\mathbf{F}}

\def\N{\mathbf{N}}
\def\Z{\mathbf{Z}}

\def\E{\mathcal{E}}

\def \l{$L_\infty$}
\def\g{\mathfrak {g}}
\def\al{\mathfrak {a}}
\def\bl{\mathfrak {b}}
\def\cl{\mathfrak {c}}
\def\dl{\mathfrak {d}}

\def\S{\mathfrak {S}}
\def\Cs{\mathfrak {C}}
\def\J{\mathfrak {J}}
\def\a{\alpha}
\def\b{\beta}
\def\V{\vartheta}
\def\w{\varpi}
\def\p{\partial}
\def\pp#1#2{\frac{\p #1}{\p #2}}
\def\PB{\left\{\cdot\,,\cdot\right\}}
\def\pb#1{\left\{#1\right\}}
\def\Pb#1{\left\{\cdot,#1\right\}}
\def\LB{\left[\cdot\,,\cdot\right]}
\def\Lb#1{\left[\cdot,#1\right]}
\def\lB#1{\left[#1,\cdot\right]}
\def\lb#1{\left[#1\right]}
\def\ideal#1{\langle #1\rangle}

\def\deg{\mathop{\rm deg}\nolimits}

\def\MC{\mathop{\mathcal MC}\nolimits}
\def\Def{\mathop{\mathcal Def}\nolimits}
\def\cla{\mathop{\rm cl}\nolimits}
\def\id{\mathop{\rm id}\nolimits}

\def\sgn{\mathop{\rm sign}\nolimits}
\def\Jac{\mathop{\rm Jac}\nolimits}
\def\Vect{{\mathfrak{X}}}

\def\L{{\mathcal L}}
\def\v{\varphi}

\def\s{\sigma}
\def\e{\varepsilon}
\def\ds{\displaystyle}
\def\Im{\mathop{\rm Im}\nolimits}

\def\Vf{{\mathcal V}}

\def\bb{{\bf b}}
\def\cc{{\bf c}}

\def\ad{\mathop{\rm ad}\nolimits}

\def\rna{\renewcommand{\arraystretch}}
\def\sc{\scriptstyle}

\newtheorem{thm}{Theorem}[section]
\newtheorem{prp}[thm]{Proposition}

\newtheorem{lma}[thm]{Lemma}
\theoremstyle{definition}

\newtheorem{rem}[thm]{Remark}

\theoremstyle{remark}

\newenvironment{eqn*}[1][1.5]
  {$$\renewcommand{\arraystretch}{#1}
      \begin{array}{rcl}}
      {\end{array}$$}

\newenvironment{eqn}[2][1.5]
  {\begin{equation}\label{#2}
   \renewcommand{\arraystretch}{#1}
   \begin{array}{rcl}}
  {\end{array}\end{equation}}

\begin{document}
\nocite{*}

\title[Deformations of Poisson structures via $L_\infty$-algebras]{$L_\infty$-interpretation of a classification of deformations of Poisson structures in dimension three}
\author{Anne Pichereau}
  \address{Max-Planck-Institut f\"ur Mathematik
       Vivatsgasse 7,
       53 111 Bonn,
       Germany}
  \email{pichereau@mpim-bonn.mpg.de}
\thanks{The author was supported by a grant of the EPDI (European Post-Doctoral Institut).}
\begin{abstract}
We give an $L_\infty$-interpretation of the classification, obtained in \cite{AP2}, of the formal deformations of a family of exact Poisson structures in dimension three. We indeed reobtain the explicit formulas for all the formal deformations of these Poisson structures, together with a classification in the generic case, by constructing a suitable quasi-isomorphism between two $L_\infty$-algebras, which are associated to these Poisson structures.
\end{abstract}
\maketitle
\tableofcontents

\section{Introduction}

In \cite{AP2}, we have exhibited a classification of the formal deformations of the Poisson structures defined on $\F[x,y,z]$ ($\F$ is an arbitrary field of characteristic zero), of the form:
\begin{equation}\label{eq:poisson}
\PB_\v = \pp{\v}{x}\,  \pp{}{y}\wedge \pp{}{z} +  \pp{\v}{y}\,  \pp{}{z}\wedge \pp{}{x}  +  \pp{\v}{z}\,  \pp{}{x}\wedge \pp{}{y} ,
\end{equation}
where $\v$ is a weight-homogeneous polynomial of $\F[x,y,z]$, admitting an isolated singularity, in the generic case.
In the present paper, following an idea of B. Fresse, we give an \l-interpretation of this result, that is to say, we obtain this result again by methods, which are different and which use the theory of \l-algebras.

The Poisson structures appear in classical mechanics, where physical systems are described by commutative algebras which are algebras of smooth functions on Poisson manifolds. They generalize the symplectic structures, as for example the natural symplectic structure on $\R^{2r}$, which was introduced by D. Poisson in 1809. On the contrary, in quantum mechanics, physical systems are described by non-commutative algebras, which are algebras of observables on Hilbert spaces, and P. Dirac has observed that, up to a factor depending on the Planck's constant, the commutator of observables appearing in the work of W. Heisenberg is the analogue of the Poisson bracket of classical mechanics. 
The Poisson structures and their deformations also appear in the theory of deformation quantization (see for instance~\cite{CKTB}) with, in particular, the very important result obtained by M.~Kontsevich in 1997:  given a Poisson manifold $(M,\pi)$ and the associative algebra $(\A=C^\infty(M),\cdot)$, there is a one-to-one correspondence between the equivalence classes of star products of~$\A$, for which the first term is $\pi$, and the equivalence classes of the formal deformations of $\pi$. 

In a more general context, a \emph{Poisson structure} on an associative commutative algebra $\A$ is a Lie algebra structure on $\A$, $\pi:\A\times\A\to \A$, which is a biderivation of $\A$ (see paragraph \ref{subsubsec:coho}).  In the case where $\A=C^\infty(M)$ is the algebra of smooth functions over a manifold $M$, one says that $(M,\pi)$ is a Poisson manifold. Formally deforming a Poisson structure $\pi$ defined on an associative commutative algebra $\A$ means considering the Poisson structures $\pi_*$ defined on the ring $\A[[\nu]]$ of all the formal power series with coefficients in $\A$ and in one parameter $\nu$, which extend the initial Poisson structure (i.e., which are $\pi$, modulo $\nu$). In this paper and in \cite{AP2}, we study a classification of formal deformations of Poisson structures modulo equivalence, two formal deformations $\pi_*$ and $\pi'_*$ of $\pi$ being equivalent if there exists a morphism $\Phi:(\A[[\nu]],\pi_*)\to(\A[[\nu]],\pi'_*)$ of Poisson algebras over $\F[[\nu]]$ which is the identity modulo $\nu$.  
There is a similar definition for the formal deformations of an associative product, the $*$-products being formal deformations of an associative product, for which each coefficient is a bidifferential operator.
We refer to \cite{AP2} for an introduction to the study of formal deformations of Poisson structures and the role played by the Poisson cohomology in this study.

\smallskip

M. Kontsevich proved the one-to-one correspondence mentioned above by using the theory of $L_\infty$-algebras and Maurer-Cartan equations. In fact, he obtained this result by proving his conjecture of \emph{formality} for a certain differential graded Lie algebra. A \emph{differential graded Lie algebra} (\emph{dg Lie algebra}, in short) is a graded Lie algebra $(\g,\LB_\g)$, endowed with a differential $\p_\g$, which is a graded derivation with respect to $\LB_\g$. The differential $\p_\g$ is a degree $1$ map satisfying $\p_\g\circ\p_\g=0$, giving rise to a cohomology $H(\g,\p_\g)$. A dg Lie algebra is a particular example of an \emph{\l-algebra} (also called \emph{strongly homotopy Lie algebra}), which is a graded vector space $L$, equipped with a collection of skew-symmetric multilinear maps $(\ell_n)_{n\in\N^*}$, satisfying different conditions, which can be viewed as generalized Jacobi identities. (The strongly homotopy algebras were introduced by J.\ Stasheff in \cite{S} in the associative case, see also \cite{LM} and \cite{LS}.)
An (\l-)\emph{quasi-isomorphism} between two dg Lie algebras (or between two \l-algebras) is an \l-morphism between them (that is to say a collection of multilinear maps $(f_n)_{n\in\N^*}$ from one to the other, satisfying a collection of compatibility conditions), which induces an isomorphism between their cohomologies. These notions will be recalled in the paragraph \ref{subsec:lalg}. 
A dg Lie algebra is said to be \emph{formal} if there exists a quasi-isomorphism between it and the dg Lie algebra given by its cohomology $H(\g,\p_\g)$ (equipped with the trivial differential and the graded Lie bracket induced by $\LB_\g$). To a dg Lie algebra $(\g,\p_\g,\LB_\g)$ is associated an equation, called the \emph{Maurer-Cartan equation} and given by:
$$
\p_\g(\gamma) + \frac{1}{2}\lb{\gamma,\gamma}_\g =0,
$$
whose solutions $\gamma\in\g^1$ are degree one homogeneous elements of $\g$, which can also be considered as depending on a formal parameter $\nu$, $\gamma\in \nu\g^1[[\nu]]$. The set of all these formal solutions is denoted by $\MC^\nu(\g)$. 
Notice that there is a also a notion of \emph{generalized Maurer-Cartan equation} associated to an \l-algebra, which is more complicated (because it takes into account the whole \l-structure).
Given a Poisson manifold $(M,\pi)$ (respectively, a Poisson algebra $(\A,\pi)$), the Poisson cohomology complex $H(M,\pi)$ associated to $(M,\pi)$ (respectively, $H(\A,\pi)$ associated to $(\A,\pi)$) is defined as follows: the cochains are the polyvector fields (respectively, the skew-symmetric multiderivations of $\A$) and the Poisson coboundary operator is given by $\delta_\pi := -\Lb{\pi}_S$, where $\LB_S$ is the Schouten bracket (which is a graded Lie bracket, obtained by extending the commutator of vector fields to a graded biderivation with respect to the wedge product). One can then associate to $\pi$ the dg Lie algebra $\g_\pi$, given by the graded vector space of the Poisson cochains (with a shift of degree), equipped with the Poisson coboundary operator associated to $\pi$  as differential (up to a sign) and the Schouten bracket as graded Lie bracket. Because a degree one element $\gamma\in\g_\pi^1$ or $\gamma\in \nu\g_\pi^1[[\nu]]$ satisfies the Jacobi identity (hence, is a Poisson structure) if and only if $\lb{\gamma,\gamma}_S=0$, an element $\gamma\in \nu\g_\pi^1[[\nu]]$ is then a formal solution of the Maurer-Cartan equation associated to~$\g_\pi$ if and only if $\pi+\gamma$ is a formal deformation of $\pi$. 
Similarly, the star products also correspond to the formal solutions of the Maurer-Cartan equation associated to a dg Lie algebra $\g_{\scriptscriptstyle H}$, constructed from the  Hochschild cohomology complex of the associative algebra $(C^\infty(M), \cdot)$. M. Kontsevich showed that the dg Lie algebra $\g_{\scriptscriptstyle H}$ is formal, by showing that it is quasi-isomorphic to the particular dg Lie algebra $\g_\pi$, associated to the trivial Poisson bracket $\pi=0$. This result, together with the fact that 
a quasi-isomorphism between two dg Lie algebras induces a bijection between the sets of all the formal solutions of the Maurer-Cartan equations modulo a gauge equivalence, leads to the desired one-to-one correspondence.

\smallskip

In this paper, we follow an idea of B. Fresse to reobtain, but with \l-methods, the explicit formulas for all the formal deformations (modulo equivalence) of $\PB_\v$ (defined in (\ref{eq:poisson}) and sometimes called \emph{exact Poisson structures}),  which were obtained in \cite{AP2}, when $\v\in\A:=\F[x,y,z]$ is a weight homogeneous polynomial with an isolated singularity, and in particular, the classification of the formal deformations of these Poisson structures, when $\v$ is generic (i.e., when its weighted degree is different from the sum of the weights of the three variables $x$, $y$ and $z$ or, equivalently, when $H^1(\A,\PB_\v)$ is zero).
To do this, we show that this classification is not a consequence of the formality of a certain dg Lie algebra, but still of the existence of a suitable \emph{quasi-isomorphism between two \l-algebras}. 
In order to explain this, let us consider $\v$ a polynomial as before and $(\g_\v,\p_\v,\LB_S)$, the dg Lie algebra associated to the Poisson algebra $(\A:=\F[x,y,z], \PB_\v)$ and $H_\v:= H(\A, \PB_\v)$ its associated Poisson cohomology. As said before, there is a shift of degree implying that $H^\ell_\v$, the homogeneous part of $H_\v$ of degree $\ell$, is in fact the $(\ell+1)$-st Poisson cohomology space $H^{\ell+1}(\A,\PB_\v)$, associated to $(\A, \PB_\v)$. 

 In fact, the classification of the formal deformations of the Poisson bracket $\PB_\v$ obtained in \cite{AP2} was indexed by elements of $H^1_\v\otimes \nu\F[[\nu]]=H^2(\A,\PB_\v)\otimes \nu\F[[\nu]]$ and B. Fresse pointed out to me that it could come from the formality of the dg Lie algebra $\g_\v$, or at least from the existence of a suitable quasi-isomorphism between $H_\v$ and the dg Lie algebra $\g_\v$, where $H_\v$ would be equipped with a suitable \l-algebra structure. In our context, having a suitable \l-algebra structure on $H_\v$ means that the set of all the formal solutions of the generalized Maurer-Cartan equation associated to $H_\v$ would be exactly $H^1_\v\otimes \nu\F[[\nu]]$. Indeed, as in the case of dg Lie algebras, a quasi-isomorphism between two \l-algebras induces an isomorphism between the sets of formal solutions of the corresponding generalized Maurer-Cartan equations, modulo a gauge equivalence, and in our case, the formal solutions of the Maurer-Cartan equation associated to $\g_\v$, modulo the gauge equivalence, correspond exactly to the equivalence classes of the formal deformations of the Poisson structure $\PB_\v$. 
 Using the explicit bases exhibited for the Poisson cohomology associated to $(\A, \PB_\v)$ in \cite{AP1} and the idea of B. Fresse, we indeed obtained the following result (see the theorems \ref{prop:defosLinfini} and \ref{prp:conclu}):
 \begin{thm}\label{prp:intro}
 Let $\v\in\A=\F[x,y,z]$ be a weight homogeneous polynomial with an isolated singularity. Consider $(\g_\v,\p_\v,\LB_S)$ the dg Lie algebra associated to the Poisson algebra $(\A,\PB_\v)$, as explained above, where $\PB_\v$ is given by (\ref{eq:poisson}), and $H_\v$ the cohomology associated to the cochain complex $(\g_\v,\p_\v)$.
 
 There exist 
 \begin{enumerate}
 \item an \l-algebra structure on $H_\v$, such that the generalized Maurer-Cartan equation associated to $H_\v$ is trivial 
 (i.e., every element $\gamma\in\nu H_\v^1[[\nu]]$ is solution);
 \item a quasi-isomorphism $f_\bullet^\v$ from the \l-algebra $H_\v$ to the dg Lie algebra $(\g_\v,\p_\v,\LB_S)$, such that the isomorphism, induced by $f_\bullet^\v$ between the formal solutions of the Maurer-Cartan equations, sends $\MC^\nu(H_\v)$ to the representatives, exhibited in \cite{AP2}, for all the formal deformations of the Poisson bracket $\PB_\v$, modulo equivalence. 
 \end{enumerate}
 \end{thm}
  This theorem \ref{prp:intro} permits us to recover the results obtained in \cite{AP2}, concerning the formal deformations of the Poisson structures $\PB_\v$.
  It also permits us to better understand  different phenomena about this result. In particular, we used in \cite{AP2} that, in the generic case, $H^1(\A,\PB_\v)$ is zero and we now know that this fact implies that the gauge equivalence in $\MC^\nu(H_\v)$ is trivial. Moreover, in the special case (when the weighted degree of $\v$ is the sum of the weights of the three variables $x$, $y$ and $z$ or, equivalently, when $H^1(\A,\PB_\v)$ is not zero), we can now better understand the equivalence classes of the formal deformations, as the equivalence relation for the formal deformations of $\PB_\v$ can be obtained by transporting the gauge equivalence in $\MC^\nu(H_\v)$ to $\MC^\nu(\g_\v)$.

\bigskip
 
 Finally, notice that, given a dg Lie algebra $(\g,\p_\g,\LB_\g)$ and a choice of bases for the cohomology spaces $H^\ell(\g,\p_\g)$ associated to the cochain complex $(\g,\p_\g)$, and using a theorem of transfer structure (see for instance the ``move" (M1) of \cite{M}), we know that there always exist an \l-algebra structure on $H(\g,\p_\g)$, together with a quasi-isomorphism between this \l-algebra and the dg Lie algebra $(\g,\p_\g,\LB_\g)$. The problem to use this result in our context where $\g=\g_\v$, is that we do not need only the existence of this \l-structure and this quasi-isomorphism, but we also need to be able to control these data, in order:
 \begin{enumerate}
 \item  for the formal solutions of the generalized Maurer-Cartan equation associated to $H_\v$ to be simple (given by $H^1_\v\otimes \nu\F[[\nu]]$),
 \item for the image of the isomorphism, induced by $f_\bullet^\v$ between the sets of formal solutions of the Maurer-Cartan equations, to give exactly the representatives of the formal deformations of the Poisson bracket $\PB_\v$ modulo equivalence which were exhibited in \cite{AP2}.
 \end{enumerate}
  To be able to do this, we have proved, in the section \ref{sec:transfer}, a proposition which permits one, given a dg Lie algebra $(\g,\p_\g,\LB_\g)$ and a choice of basis for its associated cohomology $H(\g,\p_\g)$, to construct step by step, both an \l-algebra structure $\ell_\bullet=(\ell_n)_{n\in\N^*}$ on $H(\g,\p_\g)$ with $\ell_1=0$, and a quasi-isomorphism $f_\bullet=(f_n)_{n\in\N^*}$ from $H(\g,\p_\g)$ to $\g$, such that, at each step, whatever the choices made at the previous steps for the maps $\ell_2,\dots,\ell_{m-1}$ and $f_1,\dots,f_{m-1}$, satisfying the conditions required at this step, the collections  of maps $(\ell_n)_{1\leq n\leq m-1}$ and $(f_n)_{1\leq n\leq m-1}$ extend to an \l-algebra structure $\ell_\bullet$ on $H(\g,\p_\g)$ and a quasi-isomorphism $f_\bullet$ from $H(\g,\p_\g)$ to $\g$. This result is given in the proposition \ref{prop:magic} and permits us, together with the explicit bases exhibited for the Poisson cohomology associated to $(\A, \PB_\v)$ in \cite{AP1}, to prove the desired theorem \ref{prp:intro}.

 \bigskip
 
 {\it Acknowledgments.} 
 I am grateful to B. Fresse for having introduced me to the theory of \l-algebras anf for pointing out to me a possible \l-interpretation of my results. 
 Moreover, I would like to thank J. Stasheff for interesting and useful comments on this work. 
  I also would like to thank H. Abbaspour, M. Marcolli and T. Tradler for giving me the opportunity of participating to this volume.
  Finally, this work has been done when I was a visitor at the CRM (Centre de Recerca Matematic\`a, Barcelona) and at the MPIM (Max-Planck-Institut f\" ur Mathematik, Bonn), whose hospitality are also greatly acknowledged.

\bigskip

\section{Preliminaries: $L_\infty$-algebras and Poisson algebras}

In this first section we recall some definitions in the theory of \l-algebras. We indeed need to fix our sign conventions and the notations. The notions of Poisson structures, cohomology and  the deformations of Poisson structures are also recalled.

\subsection{\l-algebras and morphisms, Maurer-Cartan equations}
\label{subsec:lalg}

We first recall the notions of \l-algebras, \l-morphisms, Maurer-Cartan equations, mainly in order to fix the sign conventions. For these notions and the conventions we choose, we refer to (the appendix A of) \cite{L} (see also \cite{LM} and \cite{DMZ}).

In this paper, $\F$ is an arbitrary field of characteristic zero and every algebra, dg Lie algebra, \l-algebra, etc, is considered over $\F$.
If $V$ is a graded\footnote{We here consider graded vector spaces as being graded over $\Z$, but for the specific cases which we study in the section \ref{sec:defosPoisson}, the considered graded vector spaces are graded only on $\N\cup\lbrace -1\rbrace$} vector space, we denote by $|x|\in \Z$ the degree of a homogeneous element $x$ of $V$. 
Let us denote by $\bigwedge^\bullet V$, the graded commutative associative algebra (rather denoted by $\bigodot^\bullet V$ in \cite{L}), obtained by dividing the tensor algebra $T^\bullet V = \bigoplus_{k\in\N} V^{\otimes k}$ of $V$ by the ideal generated by the elements of the form $x\otimes y - (-1)^{|x||y|} y\otimes x$. Denoting by $\wedge$ the product in $\bigwedge^\bullet V$, one then has:
$$
x\wedge y = (-1)^{|x||y|} y\wedge x,
$$
where $x$ and $y$ are homogeneous elements of $V$.
Then, if $x_1,\dots,x_k\in V$ are homogeneous elements of $V$ and $\s\in\S_k$ is a permutation of $\lbrace 1,\dots, k\rbrace$, one defines the so-called \emph{Koszul sign} $\e(\s;x_1,\dots,x_k)$, associated to $x_1,\dots,x_k$ and $\s$, by the equality:
$$
x_1\wedge\dots\wedge x_k = \e(\s;x_1,\dots,x_k)\; x_{\s(1)}\wedge\dots\wedge x_{\s(k)},
$$
valid in the algebra $\bigwedge^\bullet V$. Then, one also defines the number $\chi(\s;x_1,\dots,x_k)\in\lbrace -1,1\rbrace$, by:
$$
\chi(\s;x_1,\dots,x_k) := \sgn(\s)\,\e(\s;x_1,\dots,x_k),
$$
where $\sgn(\s)$ denotes the sign of the permutation $\s$. When no confusion can arise, we write $\chi(\s)$ for $\chi(\s;x_1,\dots,x_k) $. For $ i, j\in\N$, a \emph{$(i,j)$-shuffle} is a permutation $\s\in\S_{i+j}$ such that $\s(1)<\cdots <\s(i)$ and $\s(i+1)<\cdots <\s(i+j)$, and the set of all $(i,j)$-shuffles is denoted by $S_{i,j}$.\\

\subsubsection{\l-algebras}

An \emph{\l-algebra} $L$ is graded vector space $L=\bigoplus_{n\in\Z} L^n$, equipped with a collection of linear maps
$$
\ell_\bullet = \left(\ell_k : \bigotimes\nolimits^k L \to L\right)_{k\in\N^*},
$$
such that:
\begin{enumerate}
\item[a.] each map $\ell_k$ is a graded map of degree $\deg(\ell_k)=2-k$, 
\item[b.] each map $\ell_k$ is a skew-symmetric map, which means that
$$
\ell_k(\xi_{\s(1)},\dots, \xi_{\s(k)}) = \chi(\s) \ell_k(\xi_1,\dots, \xi_k),
$$ 
for all homogeneous elements $\xi_1,\dots,\xi_k\in L$ and all permutation $\s\in\S_k$;
\item[c.] the maps $\ell_k, k\in\N^*$ satisfy the following ``generalized Jacobi identity":
\begin{equation}\label{L(n)}\tag{$\J_n$}
\!\!\!\!\!\!\!\!\!\!\!\!\renewcommand{\arraystretch}{0.7} \sum_{\begin{array}{c}\scriptstyle
  i+j=n+1\\\scriptstyle i,j\geq1\end{array}}
\!\!\sum_{\s\in S_{i,n-i}}\!\! \chi(\s)(-1)^{i(j-1)}\, \ell_j\left(\ell_i\left(\xi_{\s(1)},\dots,\xi_{\s(i)}\right)\!,\xi_{\s(i+1)},\dots,\xi_{\s(n)}\right)=0,
\end{equation}
for all $n\in\N^*$ and all homogeneous $\xi_1,\dots,\xi_n\in L$.
\end{enumerate}

The map $\ell_1$ (which satisfies $\ell_1^2 =0$, by $(\J_1)$) is sometimes called the differential of $L$ and denoted by $\partial$, while the map $\ell_2$ is sometimes denoted by a bracket $\LB$.
A \emph{differential graded Lie algebra} (\emph{dg Lie algebra} in short) is an \l-algebra $(L,\ell_1,\ell_2, \ell_3,\dots)$, with $\ell_k=0$, for all $k\geq 3$.

\bigskip

Notice that if $\ell_1=0$, then the equation $(\J_n)$ reads as follows: 
$$
J_n(\L_{n-1}; \xi_1,\dots,\xi_{n})=0,
$$
 where $J_n(\L_{n-1}; \xi_1,\dots,\xi_{n})$ depends only on $\L_{n-1} := (\ell_2,\dots,\ell_{n-1})$ (and not on $\ell_n$) and is defined by:
\begin{eqnarray}\label{eq:Jn}
\lefteqn{J_n(\L_{n-1}; \xi_1,\dots,\xi_{n}):=}\nonumber\\
&&\renewcommand{\arraystretch}{0.7} \sum_{\begin{array}{c}\scriptstyle
  i+j=n+1\\\scriptstyle i,j\geq 2\end{array}}
\!\!\sum_{\s\in S_{i,n-i}} \chi(\s)(-1)^{i(j-1)}\, \ell_j\left(\ell_i\left(\xi_{\s(1)},\dots,\xi_{\s(i)}\right)\!,\xi_{\s(i+1)},\dots,\xi_{\s(n)}\right).
\end{eqnarray}
When no confusion can arise, we rather write $J_n(\xi_1,\dots,\xi_{n})$ for $J_n(\L_{n-1}; \xi_1,\dots,\xi_{n})$.
\\

\subsubsection{\l-morphisms, quasi-isomorphisms}

There is a notion of (weak) morphism of \l-algebras, which we do not need here. We only need the particular case when the considered morphism goes from an \l-algebra to a dg Lie algebra. 
(For the general definition of \l-morphisms between \l-algebras, see \cite{KS}.)
Let $L=(L,\ell_1,\ell_2,\dots)$ be an \l-algebra and let $\g=(\g,\p_\g,\LB_\g)$ be a dg Lie algebra. A \emph{(weak) \l-morphism} from $L$ to $\g$ is a collection of linear maps 
$$
f_\bullet = \left(f_n : \bigotimes\nolimits^n L \to \g\right)_{n\in\N^*},
$$
such that:
\begin{enumerate}
\item[a.] each map $f_n$ is a graded map of degree $\deg(f_n)=1-n$;
\item[b.] each map $f_n$ is skew-symmetric;
\item[c.] the following identities hold, for all $n\in\N^*$ and all homogeneous elements $\xi_1,\dots,\xi_n\in L$:
\begin{equation}\label{E(n)}
   \renewcommand{\arraystretch}{1.7}
     \begin{array}{l}
\quad\p_\g\left(f_n(\xi_1,\dots,\xi_n)\right)\\
+ \!\!\!
\ds\sum_{\renewcommand{\arraystretch}{0.5}
     \begin{array}{c}
          \scriptstyle j+k=n+1\\
          \scriptstyle j,k\geq1
     \end{array}}
 \sum_{\s\in S_{k,n-k}} \chi(\s)\,(-1)^{k(j-1)+1}\, f_j\left(\ell_k\left(\xi_{\s(1)},\dots,\xi_{\s(k)}\right)\!,\xi_{\s(k+1)},\dots,\xi_{\s(n)}\right)\\
+ \!\!\!
\ds\sum_{\renewcommand{\arraystretch}{0.5}
    \begin{array}{c}
        \scriptstyle s+t=n\\
        \scriptstyle s,t\geq1\\
     \end{array}}
\!\!\!\!\sum_{\renewcommand{\arraystretch}{0.5}
    \begin{array}{c}
        \scriptstyle \tau\in S_{s,n-s}\\
        \scriptstyle \tau(1)<\tau(s+1)
     \end{array}} \!\!\! \!\!\! \!
  \chi(\tau) \,e_{s,t}(\tau)
      \lb{f_s \left(\xi_{\tau(1)},\dots,\xi_{\tau(s)}\right), f_t \left(\xi_{\tau(s+1)},\dots,\xi_{\tau(n)}\right)}_\g =0,
     \end{array}
\end{equation}
where $e_{s,t}(\tau) :=(-1)^{s-1}\cdot(-1)^{(t-1)\left(\sum\limits_{p=1}^s |\xi_{\tau(p)}|\right)}$ and where $\chi(\s)$ (respectively, $\chi(\tau)$) stands for $\chi(\s;\xi_1,\dots,\xi_n)$ (respectively, $\chi(\tau;\xi_1,\dots,\xi_n)$).
\end{enumerate}

\smallskip

We point out that, for $1\leq s\leq n$ and for an $(s,n-s)$-shuffle $\tau\in S_{s,n-s}$, the condition $ \tau(1)<\tau(s+1)$ is equivalent to $\tau(1) =1$.
One says that the \l-morphism $f_\bullet$ from $L$ to $\g$ is a \emph{quasi-isomorphism} (or a \emph{(weak) \l-equivalence}) if the chain map 
$f_1: (L, \ell_1)\to (\g,\p_\g)$ induces an isomorphism between the cohomologies associated to the cochain complexes $(L, \ell_1)$ and $(\g,\p_\g)$.

\bigskip

Notice that if the \l-algebra $(L,\ell_1,\ell_2,\dots)$ satisfies $\ell_1=0$, then we write the equation (\ref{E(n)}) rather in the following form:
\begin{equation}\label{E'(n)}\tag{$\E_n$}
\p_\g\left(f_n(\xi_1,\dots,\xi_n)\right) - f_1\left(\ell_n\left(\xi_1,\dots,\xi_n\right)\right) = T_n(\Fe_n, \L_{n-1}; \xi_1,\dots,\xi_n)
\end{equation}
where $T_n(\Fe_n, \L_{n-1}; \xi_1,\dots,\xi_n)$ depends on the elements $\Fe_n:=(f_1,\dots,f_{n-1})$ and $\L_{n-1}:=(\ell_2,\dots,\ell_{n-1})$, and is defined by:
\begin{equation}\label{T(n)}
   \renewcommand{\arraystretch}{1.7}
     \begin{array}{l}
\quad T_n(\Fe_n, \L_{n-1}; \xi_1,\dots,\xi_n) :=\\
\!\!\!
\ds\sum_{\renewcommand{\arraystretch}{0.5}
     \begin{array}{c}
          \scriptstyle j+k=n+1\\
          \scriptstyle j,k\geq2
     \end{array}}
 \sum_{\s\in S_{k,n-k}} \chi(\s)\,(-1)^{k(j-1)}\, f_j\left(\ell_k\left(\xi_{\s(1)},\dots,\xi_{\s(k)}\right)\!,\xi_{\s(k+1)},\dots,\xi_{\s(n)}\right)\\
- \!\!\!
\ds\sum_{\renewcommand{\arraystretch}{0.5}
    \begin{array}{c}
        \scriptstyle s+t=n\\
        \scriptstyle s,t\geq1\\
     \end{array}}
\!\!\!\!\sum_{\renewcommand{\arraystretch}{0.5}
    \begin{array}{c}
        \scriptstyle \tau\in S_{s,n-s}\\
        \scriptstyle \tau(1)=1
     \end{array}} \!\!\! \!\!\! \!
  \chi(\tau) \,e_{s,t}(\tau)
      \lb{f_s \left(\xi_{\tau(1)},\dots,\xi_{\tau(s)}\right), f_t \left(\xi_{\tau(s+1)},\dots,\xi_{\tau(n)}\right)}_\g,
     \end{array}
\end{equation}
 for all $n\in\N^*$ and all (homogeneous) elements $\xi_1,\dots,\xi_n\in L$. When no confusion can arise, we simply write $T_n(\xi_1,\dots,\xi_n)$ for $T_n(\Fe_n, \L_{n-1}; \xi_1,\dots,\xi_n)$.
  \\

\subsubsection{Maurer-Cartan equation} 

To an \l-algebra is associated the so-called generalized Maurer-Cartan equation (or homotopy Maurer-Cartan equation). In our context, we only need a particular case of it, where the solutions depend formally on a parameter $\nu$. 
Let $L=(L,\ell_1,\ell_2,\dots)$ be an \l-algebra. The \emph{generalized Maurer-Cartan equation} associated to $L$ is written as follows:
\begin{equation}\label{MC}
-\ell_1(\gamma) - \frac{1}{2} \ell_2(\gamma,\gamma) + \frac{1}{3!} \ell_3(\gamma,\gamma,\gamma) + \cdots + \frac{(-1)^{n(n+1)/2}}{n!} \ell_n(\gamma,\dots,\gamma) + \cdots  =0,
\end{equation}
for $\gamma \in L^1\otimes \nu\F[[\nu]] = \nu L^1[[\nu]]$, where $\nu$ is a formal parameter. Notice that the maps $\ell_n$, $n\in\N^*$ are extended by multilinearity with respect to the parameter $\nu$ (and are still denoted by $\ell_n$) and that this infinite sum (\ref{MC}) is well-defined because there is no constant term in $\gamma$ (i.e., $\gamma$ is zero modulo $\nu$), so that the coefficient of each $\nu^i$, $i\in\N$ is given by a finite sum. 
The same will hold for the equations (\ref{gauge}) and (\ref{eq:MCf}).

The set of all the solutions of the generalized Maurer-Cartan equation associated to $L$ and depending formally on a parameter~$\nu$ is denoted by $\MC^\nu(L)$.
One introduces the \emph{gauge equivalence} on this set, which is denoted by $\sim$ and generated by infinitesimal transformations of the form:
\begin{equation}\label{gauge}
\gamma \longmapsto \xi\cdot\gamma := \gamma - \sum_{n\in\N^*}  \frac{(-1)^{n(n-1)/2}}{(n-1)!} \ell_n(\xi,\gamma,\gamma,\dots,\gamma),
\end{equation}
where $\xi\in L^0\otimes \nu\F[[\nu]]$.
\begin{rem}\label{rem:MCg}
Let us consider the particular case where the \l-algebra $L$ is a dg Lie algebra $(\g,\p_\g,\LB_\g)$ whose differential is given by $\p_\g = \lB{\chi}_\g$, for some degree one element $\chi\in\g^1$ satisfying $\lb{\chi,\chi}_\g=0$. Then, we have:
\begin{eqn}{MCg}
\MC^\nu(\g) &=& \left\lbrace \gamma \in \g^1\otimes \nu\F[[\nu]] \mid 
                 \lb{\chi,\gamma}_\g + \frac{1}{2} \lb{\gamma,\gamma}_\g=0\right\rbrace\\
&=& \left\lbrace \gamma \in \g^1\otimes \nu\F[[\nu]] \mid 
                 \lb{\chi + \gamma, \chi + \gamma}_\g=0\right\rbrace.
\end{eqn}
Moreover the infinitesimal transformation (\ref{gauge}) becomes in this case:
\begin{equation}\label{gauge-g}
\gamma \longmapsto 
\xi\cdot\gamma := \gamma  + \lb{\xi,\chi + \gamma}_\g,
\end{equation}
for $\xi\in \g^0\otimes \nu\F[[\nu]]$.
\end{rem}

We denote by $\Def^\nu(L)$ the set of all the gauge equivalence classes of the formal solutions of the generalized Maurer-Cartan equation associated to $L$, 
$$
\Def^\nu(L) := \MC^\nu(L) / \sim.
$$
 For $\gamma\in  \MC^\nu(L)$, we denote by $\cla(\gamma)\in\Def^\nu(L)$ its equivalence class modulo the gauge equivalence.
In the following we will use the theorem (see for instance \cite{K} or~\cite{DMZ}):
\begin{thm}[\cite{K}, \cite{DMZ}, \dots]\label{thmMC}
Let $L$ and $L'$ be two \l-algebras and let us suppose that $f_\bullet = \left(f_n : \bigotimes\nolimits^n L \to L'\right)_{n\in\N^*}$
is a quasi-isomorphism from $L$ to $L'$. Then $f_\bullet$ induces an isomorphism $\Def^\nu(f_\bullet)$ from $\Def^\nu(L)$ to $\Def^\nu(L')$. This isomorphism is given, for $\gamma\in \MC^\nu(L)$, by:
\begin{equation}\label{eq:def(f)}
\Def^\nu(f_\bullet)\left(\cla(\gamma)\right) := \cla\left(\MC^\nu(f_\bullet)(\gamma) \right),
\end{equation}
where 
\begin{equation}\label{eq:MCf}
\MC^\nu(f_\bullet)(\gamma) :=\sum\limits_{n\geq 1} \frac{(-1)^{1+n(n+1)/2}}{n!} f_n(\gamma,\dots,\gamma).
\end{equation}
\end{thm}
Notice that we will only use this theorem in the case $L'$ is a dg Lie algebra.

\subsection{Poisson algebras, cohomology and deformations}

In this paper, our goal is to apply the theory of \l-algebras to the problem of deformations of Poisson structures. We here recall the notions of Poisson algebras, cohomology and deformations, and explain how one can associate a dg Lie algebra to a Poisson algebra.

\subsubsection{Poisson algebra and cohomology}
\label{subsubsec:coho}

We recall that a \emph{Poisson
structure} $\PB$ (also denoted by $\pi_0$) on an associative
commutative algebra $(\A,\cdot)$ is a skew-symmetric biderivation of
$\A$, i.e., a map $\PB : \bigwedge^2\A\to\A$ satisfying the derivation
property:
\begin{equation}\label{leibniz}
\pb{FG,H} =F\pb{G,H} + G\pb{F,H}, \hbox{ for all } F,G,H\in\A,
\end{equation}
(where $FG$ stands for $F\cdot G$), 
which is also a Lie structure on $\A$, i.e., which satisfies the Jacobi identity
\begin{equation}\label{leibniz}
\pb{\pb{F,G},H}+ \pb{\pb{G,H},F} + \pb{\pb{H,F}, G}=0, \hbox{ for all } F,G,H\in\A.
\end{equation}
 The couple $(\A,\PB=\pi_0)$ is then called a \emph{Poisson algebra}.

The Poisson cohomology has been introduced by A.\ Lichnerowicz in \cite{Li}; see also \cite{H2} for an algebraic approach.
The \emph{Poisson cohomology} complex, associated to a Poisson algebra $(\A,\pi_0)$, is defined as follows.
The space of all Poisson
cochains is $\Vect^\bullet(\A):=\bigoplus_{k\in\N}\Vect^k(\A)$, where
$\Vect^0(\A)$ is $\A$ and, for all $k\in\N^*$, $\Vect^k(\A)$ denotes
the space of all skew-symmetric $k$-derivations of $\A$, i.e., the skew-symmetric
$k$-linear maps $\A^k\to\A$ that satisfy the derivation property
(\ref{leibniz}) in each of their arguments.
The Poisson coboundary operator
$\delta^k_{\pi_0}:\Vect^k(\A)\to\Vect^{k+1}(\A)$ is given by the formula
$$
\delta^k_{\pi_0} := -\Lb{\pi_0}_S,
$$ 
where $\LB_S:\Vect^p(\A)\times\Vect^q(\A)\to\Vect^{p+q-1}(\A)$ is the so-called
Schouten bracket. The Schouten bracket is a graded
Lie bracket that generalizes the commutator of derivations and that is
a graded biderivation with respect to the wedge product of multiderivations (see \cite{PLV}). It is
defined, for $P\in\Vect^p(\A)$, $Q\in\Vect^q(\A)$ and for
$F_1,\dots,F_{p+q-1}\in\A$, by:
\begin{eqnarray}\label{schou_exp}
\lefteqn{\lb{P,Q}_S[F_1,\dots,F_{p+q-1}]}\hspace{0cm} \nonumber\\
    &&\quad= \sum_{\s\in S_{q,p-1}}\!\!\sgn(\s)P\left[Q[F_{\s(1)},\dots,F_{\s(q)}],
                  F_{\s(q+1)},\dots,F_{\s(q+p-1)}\right]\\
    &&- (-1)^{(p-1)(q-1)}\sum_{\s\in S_{p,q-1}}\!\!\sgn(\s)Q\left[P[F_{\s(1)},\dots,F_{\s(p)}],
                  F_{\s(p+1)},\dots,F_{\s(p+q-1)}\right].\nonumber
\end{eqnarray}
It is easy and useful to verify that, given a skew-symmetric biderivation $\pi\in\Vect^2(\A)$, the Jacobi identity for $\pi$ is equivalent to $\lb{\pi, \pi}_S=0$, in other words, if  $\pi\in\Vect^2(\A)$ is a skew-symmetric biderivation of $\A$, then 
$\pi$ is a Poisson structure on $\A$ if and only if  $\lb{\pi,\pi}_S=0$.
\\

\subsubsection{The dg Lie algebra associated to the Poisson complex}
\label{sbsec:dglie-poisson}

The Poisson cohomology complex associated to a Poisson algebra $(\A,\pi_0)$ together with the Schouten bracket give rise to a dg Lie algebra, $(\g,\p_\g,\LB_\g)$, defined as follows.
\begin{enumerate}
\item For all $n\in\N^*$, the degree $n$ homogeneous part of $\g$ is given by 
$$
\g^n := \Vect^{n+1}(\A),
$$
so that the degree of $P\in\Vect^p(\A)=\g^{p-1}$, viewed as an element of $\g$, is $|P|:=p-1$,
\item for all $P\in\Vect^p(\A)=\g^{p-1}$, 
$$
\p_\g (P) := (-1)^{|P|} \delta_{\pi_0}^p(P) = (-1)^{p-1} \delta_{\pi_0}^p(P),
$$
\item the graded Lie bracket on $\g$ is given by the Schouten bracket:
$$
\LB_\g := \LB_S.
$$
\end{enumerate}
Notice that, using the skew-symmetry of the Schouten bracket, and the definition of $\delta_{\pi_0}^p$, we can write $\p_\g = \lB{\pi_0}_S$.
As $\lb{\pi_0,\pi_0}_S=0$ (see last paragraph \ref{subsubsec:coho}), the dg Lie algebra $(\g,\p_\g,\LB_\g)$ associated to a Poisson algebra $(\A,\pi_0)$ satisfied the conditions of the remark \ref{rem:MCg}.
\\

\subsubsection{Formal deformations of Poisson structures}
\label{par:formal-defos}

In this paragraph, we define the notion of formal deformations of Poisson structures. For more details about this, see \cite{AP2}.
Let $(\A,\cdot)$  be an associative commutative algebra over $\F$ and let $\pi_0$ be a Poisson structure on $(\A,\cdot)$.
We consider the $\F[[\nu]]$-vector
space $\A[[\nu]]$ of
all formal power series in $\nu$, with coefficients in $\A$. The
associative commutative product~``$\cdot$'', defined on $\A$, is naturally extended to
an associative commutative product on $\A[[\nu]]$, still denoted by
``$\cdot$''. 
A \emph{formal deformation} of $\pi_0$ is a Poisson structure on the associative
$\F[[\nu]]$-algebra $\A[[\nu]]$, that extends the initial
Poisson structure. In other words, it is given by a map 
$\pi_* : \A[[\nu]]\times\A[[\nu]] \to \A[[\nu]]$ satisfying the Jacobi
identity and of the form:
$$
\pi_* = \pi_0 +\pi_1\nu +\cdots +\pi_n\nu^n +\cdots,
$$
where the $\pi_i$ are skew-symmetric biderivations of $\A$ (extended by bilinearity with respect to $\nu$).
Notice that given a map $\pi_* = \pi_0 +\pi_1\nu +\cdots +\pi_n\nu^n
+\cdots : \A[[\nu]]\times\A[[\nu]]\to\A[[\nu]]$ where for all $i\in\N$,
$\pi_i\in\Vect^2(\A)$ is a skew-symmetric biderivation of $\A$, we
have that $\pi_*$ is a formal deformation of $\pi_0$ if and only if
$\lb{\pi_*,\pi_*}_S=0$.

There is a natural notion of equivalence for deformations of a Poisson structure~$\pi_0$. 
Two formal deformations $\pi_*$ and $\pi'_*$ of $\pi_0$ are said to be
\emph{equivalent} if there exists
an $\F[[\nu]]$-linear map $\Phi:(\A[[\nu]],\pi_*)\to(\A[[\nu]],\pi_*')$, which is equal to the identity modulo $\nu$ and is a Poisson morphism, i.e., it is a morphism of associative algebras $\Phi:(\A[[\nu]],\cdot)\to(\A[[\nu]],\cdot)$, which satisfies:
\begin{equation}\label{eq:Phi}
\pi'_*[\Phi(F),\Phi(G)] = \Phi(\pi_*[F,G]),
\end{equation}
for all $F,G\in\A$ (and therefore, for all $F,G\in\A[[\nu]]$). It is also possible to write such a morphism $\Phi$ as the exponential of an element $\xi\in \nu\Vect^1(\A)[[\nu]] = \Vect^1(\A)\otimes \nu\F[[\nu]]$, so that  (see for example the lemma 2.1 of \cite{AP2}) the map $\pi'_*$ given by (\ref{eq:Phi}) can also be written as:
\begin{equation}\label{eq:exp(adxi)}
\pi'_* = e^{\ad_\xi}(\pi_*) 
= \ds \pi_* + \sum_{k\in\N^*} \frac{1}{k!} 
\underbrace{\lb{\xi,\lb{\xi,\dots\,,\lb{\xi,\pi_*}_S\dots}_S}_S}
_{\hbox{$k$ brackets}}.
\end{equation}%

Let us now consider the dg Lie algebra $(\g,\p_\g,\LB_\g)$ associated to the Poisson algebra $(\A,\pi_0)$, as explained in the previous paragraph \ref{sbsec:dglie-poisson}.
According to the remark \ref{rem:MCg}, we have:
\begin{eqn*}
\MC^\nu(\g)  &=& 
             \left\lbrace \gamma = \sum\limits_{i\geq 1}\pi_i\nu^i\in \Vect^2(\A)\otimes \nu\F[[\nu]] \mid\right.\\
                           & & \qquad    \left.   \pi_*:=\pi_0+\sum\limits_{i\geq 1}\pi_i\nu^i \hbox{ is a formal deformation of } \pi_0\right\rbrace,
\end{eqn*}%
so that, there is a natural one-to-one correspondence between $\MC^\nu(\g)$ and the space of all formal deformations of $\pi_0$:
$$
\begin{array}{rcl}
\MC^\nu(\g) &\to& \left\lbrace \hbox{formal deformations of } \pi_0\right\rbrace\\
\gamma &\mapsto& \pi_0+\gamma
\end{array}
$$
Moreover, the infinitesimal transformation (\ref{gauge-g}) on elements of $\MC^\nu(\g)$ can be transposed to an infinitesimal transformation on formal deformations $\pi_*$ of $\pi_0$. It then becomes:
$$
\pi_* \;\longmapsto\; \xi\cdot \pi_* := \pi_*+ \lb{\xi,\pi_*}_S,
$$
for $\xi\in \Vect^1(\A)\otimes \nu\F[[\nu]]$. We conclude that there is a one-to-one correspondence between the elements of $\Def^\nu(\g)$ and the equivalence classes of the formal deformations of $\pi_0$.

\section{Choice in a transfer of \l-algebra structure}
\label{sec:transfer}

For $\g =(\g,\p_\g,\LB_\g)$ a dg Lie algebra, we denote by $H(\g,\p_\g)$, the graded vector space given by the cohomology of the cochain complex $(\g,\p_\g)$. Equipped with the trivial differential, it is a cochain complex $(H(\g,\p_\g),0)$. Moreover, $Z(\g,\p_\g)$ denotes the graded vector space of all the cocycles of the cochain complex $(\g,\p_\g)$: 
$$
Z(\g,\p_\g) = \ker \p_\g \subseteq \g,
$$
 and $B(\g,\p_\g)$, the graded vector space of all its coboundaries: 
$$
B(\g,\p_\g) = \Im\p_\g \subseteq \g,
$$
so that $H(\g,\p_\g) = Z(\g,\p_\g) / B(\g,\p_\g)$. (The grading of $Z(\g,\p_\g)$, $B(\g,\p_\g)$ and $H(\g,\p_\g)$ is naturally induced by the grading of $\g$.)
We denote by $p$ the natural projection from $Z(\g,\p_\g)$ to the cohomology of $\g$, and for every cocycle $x\in Z(\g,\p_\g)\subseteq\g$, the notations $p(x)$ and $\bar x$ both stand for the cohomological class of $x$,
\begin{equation}\label{eqp}
  \begin{array}{lcccl}
   p &:& Z(\g,\p_\g)&\to&H(\g,\p_\g)\\
     & & x &\mapsto& p(x) =\bar x.
  \end{array}
\end{equation}
We now define a graded linear map $f_1$, of degree $0$, from $H(\g,\p_\g)$ to $\g$. This definition depends on a choice of a basis $\bb^\ell$, for each cohomology space $H^\ell(\g,\p_\g)$, and on a choice of representatives $\left(\V^\ell_k\right)_k$ of the elements of the basis $\bb^\ell$:
$$
\bb^\ell =\left(\overline{\V^\ell_k}\right)_k.
$$
(We do not need here to specify the set by which the basis $\bb^\ell$ is indexed.)
Then the map $f_1: H(\g,\p_\g) \to \g$ is defined by
\begin{equation}\label{eqf_1}
\renewcommand{\arraystretch}{1.5}
  \begin{array}{lcccl}
    f_1 &:& H^\ell(\g,\p_\g)&\to&Z^\ell(\g,\p_\g) \subseteq\g^\ell\\
     & & \xi = \sum_{k}\lambda^\ell_k\, \overline{\V^\ell_k}  &\mapsto& \sum_{k}\lambda^\ell_k\, \V^\ell_k,
  \end{array}
 \end{equation}
for all $\ell\in\Z$, and where $\xi = \sum_{k}\lambda^\ell_k\, \overline{\V^\ell_k}$ is the unique decomposition of 
$\xi\in H^\ell(\g,\p_\g)$ in the fixed basis $\bb^\ell$ (the $\lambda^\ell_k$ are constants). We deduce from the definition of $f_1$ that we have:
\begin{equation}\label{eq:decomp-f1}
Z(\g,\p_\g) \simeq \Im f_1 \oplus B(\g,\p_\g),
\end{equation}
and 
\begin{equation}\label{eq:property-f1}
x - f_1\circ p (x) \in  B(\g,\p_\g), \quad \hbox{ for all } x\in Z(\g,\p_\g).
\end{equation}
Also the map $f_1$ is a chain map between the two cochain complexes $(H(\g,\p_\g), 0)$ and $(\g,\p_\g)$, which induces an isomorphism between their cohomologies. This implies in particular that if one extends $f_1$ to a (weak) \l-morphism 
$$
f_\bullet =\left(f_n: \bigotimes\nolimits^n H(\g,\p_\g) \to \g\right)_{n\in\N^*}
$$
 (where $H(\g,\p_\g)$ is equipped with an \l-algebra structure), then $f_\bullet$ is automatically a quasi-isomorphism.

We indeed want to construct an \l-algebra structure on $H(\g,\p_\g)$ together with a quasi-isomorphism from it to the dg Lie algebra $\g$. We know that, by using a theorem of \l-algebra structure transfer, (see for instance the ``move" (M1) of \cite{M}), there exists such a \l-algebra structure on $H(\g,\p_\g)$ and such a quasi-isomorphism from it to the dg Lie algebra $\g$, which extends $f_1$, but, as explained in the introduction, the point here is that we need to construct a \emph{specific} \l-algebra structure on $H(\g,\p_\g)$ and a \emph{specific} quasi-isomorphism. This prevents one to express the transfer structure in terms of a homotopy map (as usually done with the pertubation lemma) because it seems to the author that such a map cannot be explicitly written in general and especially in the context we will use in the section \ref{sec:defosPoisson}.
In order to have as much control in this contruction as possible, we show  the following:
\begin{prp}\label{prop:magic}
Let $\g =(\g,\p_\g,\LB_\g)$ be a dg Lie algebra, let $H(\g,\p_\g)$ denote the graded space given by the cohomology associated to the cochain complex $(\g,\p_\g)$.
We fix $f_1 : H(\g,\p_\g)\to\g$ as being the map defined in (\ref{eqf_1}), associated to a choice of bases $(\bb^\ell)_\ell$ for the cohomology spaces $(H^\ell(\g,\p_\g))_\ell$.
We also fix $\ell_1 : H(\g,\p_\g)\to H(\g,\p_\g)$ as being trivial ($\ell_1=0$) so that the equations $(\E_1)$ and $(\J_1)$ are automatically satisfied.
\begin{enumerate}
\item[(a)] There exist skew-symmetric graded linear maps
\begin{equation*}
\ell_2 :  H(\g,\p_\g)\otimes H(\g,\p_\g) \to H(\g,\p_\g), \quad\hbox{ and }\quad 
f_2 : H(\g,\p_\g)\otimes H(\g,\p_\g) \to \g,
\end{equation*}
of degrees $\deg(\ell_2)=0$ and $\deg(f_2)=-1$, such that the equations $(\E_2)$ and $(\J_2)$ are satisfied. Moreover, such a map $\ell_2$ satisfies also the equation~$(\J_3)$.
\item[(b)] Let $m\geq 3$ be an integer. For any skew-symmetric graded linear maps 
\begin{eqnarray*}
\ell_k &:& \bigotimes\nolimits^k H(\g,\p_\g) \to H(\g,\p_\g), \qquad \hbox{ for } 2\leq k\leq m-1,\\
f_k &: &\bigotimes\nolimits^k H(\g,\p_\g) \to \g, \qquad \hbox{ for } 2\leq k\leq m-1,
\end{eqnarray*}
of degrees $\deg(\ell_k)=2-k$ and $\deg(f_k)=1-k$, for all $2\leq k\leq m-1$, and such that the equations $(\J_2)$ -- $(\J_m)$ and $(\E_2)$ -- $(\E_{m-1})$ are satisfied, there exist skew-symmetric graded linear maps 
$$
\ell_m: \bigotimes\nolimits^m H(\g,\p_\g) \to H(\g,\p_\g) \quad\hbox{ and }\quad f_m : \bigotimes\nolimits^m H(\g,\p_\g) \to \g,
$$ 
with $\deg(f_m) = 1-m$, $\deg(\ell_m) = 2-m$ and satisfying the equation~$(\E_{m})$. Moreover, such a map $\ell_m$ necessarily satisfies also the equation $(\J_{m+1})$.
\end{enumerate}
\end{prp}
\begin{rem}
This proposition implies in particular that there exist an \l-algebra structure $\ell_\bullet$ on $H(\g,\p_\g)$ with the trivial differential $\ell_1=0$ and a quasi-isomorphism $f_\bullet $ from $H(\g,\p_\g)$ to $\g$ that extends $f_1$ (defined in (\ref{eqf_1})). But, this proposition implies morever that, whatever the choices made for the first $m-1$ maps $\ell_1,\dots,\ell_{m-1}$ and $f_1,\dots, f_{m-1}$ ($m$ is an arbitrary integer), with $\ell_1=0$ and $f_1$ given by (\ref{eqf_1}), if these maps satisfy the first $m$ equations defining an \l-algebra structure (equations $(\J_1)$ -- $(\J_m)$) and the first $m-1$ equations defining an \l-morphism (equations $(\E_1)$ -- $(\E_{m-1})$), then they still extend to an \l-algebra structure $\ell_\bullet$ on $H(\g,\p_\g)$ and a quasi-isomorphism $f_\bullet $ from $H(\g,\p_\g)$ to $\g$.
\end{rem}
\begin{proof}
The idea of this proof is similar to the one used by T.\ Kadeishvili in \cite{Ka} (where he considers $A_\infty$-algebras) to prove his theorem 1.
Let us first prove the part (a) of this proposition. To do this, we first show (\textit{Step 1}) that the identity $\ell_1=0$ and the definition (\ref{eqf_1}) of $f_1$ imply that $\p_\g\left(T_2(\xi_1,\xi_2)\right)=0$, for all $\xi_1,\xi_2\in H(\g,\p_\g)$. By (\ref{eq:decomp-f1}), the cocycle $T_2(\xi_1,\xi_2)$ then decomposes  as a coboundary (element in the image of $\p_\g$) plus an element in the image of $f_1$, which permit us to conclude the existence of both maps $f_2$ and $\ell_2$, satisfying the equation ($\E_{2}$). Secondly (\textit{Step 2}), we show that the obtained map $\ell_2$, satisfying $(\E_2)$, also necessarily satisfies $(\J_3)$. 

\bigskip

(a) - \textit{ \underline{Step 1.}}
The skew-symmetric graded linear maps $f_1$ (given by (\ref{eqf_1})) and $\ell_1:=0$ are of degree $0$ and $-1$ respectively, and satisfy both equations:
\begin{equation}\label{J'(1)}\tag{$\J_1$}
\ell_1 \circ \ell_1 =0,   
\end{equation}
and
\begin{equation}\label{E'(1)}\tag{$\E_1$}
\p_\g \circ f_1 =0.
\end{equation}
Let $\xi_1,\xi_2\in H(\g,\p_\g)$. We have $T_2(\xi_1,\xi_2)= -\lb{f_1(\xi_1),f_1(\xi_2)}_\g$. As $(\g,\p_\g,\LB_\g)$ is a dg Lie algebra, $\p_\g$ is a (graded) derivation for $\LB_\g$, hence:
 \begin{equation}\label{eq:dT2}
 \p_\g\left(T_2(\xi_1,\xi_2)\right)  = - \lb{ \p_\g\left(f_1(\xi_1)\right) ,f_1(\xi_2)}_\g - (-1)^{|\xi_1|}\lb{f_1(\xi_1),\p_\g\left(f_1(\xi_2)\right)}_\g
      = 0,
\end{equation}
by (\ref{E'(1)}). We now define  a skew-symmetric graded linear map $\ell_2 : \bigwedge\nolimits^2 H(\g,\p_\g)\to H(\g,\p_\g)$ of degree $0$,  by:
\begin{equation}\label{eq:ell2}
\ell_2(\xi_1,\xi_2) :=-p\circ T_2(\xi_1,\xi_2),
\end{equation}
for all $\xi_1,\xi_2\in H(\g,\p_\g)$. This map is well-defined because, according to (\ref{eq:dT2}), $T_2(\xi_1,\xi_2)$ is a cocycle for the cochain complex $(\g,\p_\g)$, and it trivially satisfies the equation $(\J_2)$, because $\ell_1=0$. It is also possible, according to (\ref{eq:property-f1}), to define a skew-symmetric graded linear map $f_2:\bigwedge\nolimits^2 H(\g,\p_\g)\to \g$, of degree $-1$, with the following formula:
\begin{equation}
\p_\g\left(f_2(\xi_1,\xi_2)\right) = T_2(\xi_1,\xi_2) - f_1\circ p\left(T_2(\xi_1,\xi_2)\right),
\end{equation}
for all $\xi_1,\xi_2\in H(\g,\p_\g)$. The maps $\ell_2$ and $f_2$ then satisfy the equation $(\E_2)$, because $-f_1\circ p\left(T_2(\xi_1,\xi_2)\right) = f_1\circ\ell_2(\xi_1,\xi_2)$.
 Notice that, for every $\xi_1,\xi_2\in H(\g,\p_\g)$, the choice of the element $f_2(\xi_1,\xi_2)\in\g$ is unique, up to a cocycle.

\bigskip

(a) - \textit{ \underline{Step 2.}} Now, let us prove the second part of (a), by showing that the map $\ell_2$, defined in (\ref{eq:ell2}), satisfies the equation
\begin{equation}\label{L(3)}\tag{$\J_3$}
\renewcommand{\arraystretch}{0.7}\sum_{\s\in S_{2,1}} \chi(\s)\, \ell_2\left(\ell_2\left(\xi_{\s(1)},\xi_{\s(2)}\right)\!,\xi_{\s(3)}\right)=0,
\end{equation}
for all homogeneous $\xi_1,\xi_2,\xi_3\in H(\g,\p_\g)$, where $\chi(\s)$ stands for $\chi(\s;\xi_1,\xi_2,\xi_3)$.
We prove this, by using the equations ($\E_1$) and ($\E_2$) and the graded Jacobi identity satisfied by $\LB_\g$.
Let $\xi_1,\xi_2,\xi_3\in H(\g,\p_\g)$ and let $\s\in S_{2,1}$. By the definition (\ref{eq:ell2}) of $\ell_2$, we have 
$\ell_2\left(\ell_2\left(\xi_{\s(1)},\xi_{\s(2)}\right)\!,\xi_{\s(3)}\right) = -p\left(T_2\left(\ell_2\left(\xi_{\s(1)},\xi_{\s(2)}\right)\!,\xi_{\s(3)}\right)\right)$. Moreover, by definition of $T_2$, 
\begin{eqnarray*}
\lefteqn{T_2\left(\ell_2\left(\xi_{\s(1)},\xi_{\s(2)}\right)\!,\xi_{\s(3)}\right) = -\lb{f_1\left(\ell_2\left(\xi_{\s(1)},\xi_{\s(2)}\right)\right)\!,f_1\left(\xi_{\s(3)}\right)}_\g}\\
 &=& \lb{T_2 \left(\xi_{\s(1)},\xi_{\s(2)} \right),f_1\left(\xi_{\s(3)}\right)}_\g -  \lb{\p_\g\left(f_2 \left(\xi_{\s(1)},\xi_{\s(2)} \right)\right)\!,f_1\left(\xi_{\s(3)}\right)}_\g\\
 &=&  -\lb{\lb{f_1\left(\xi_{\s(1)}\right),f_1\left(\xi_{\s(2)}\right)}_\g,f_1\left(\xi_{\s(3)}\right)}_\g -  \lb{\p_\g\left(f_2 \left(\xi_{\s(1)},\xi_{\s(2)} \right)\right)\!,f_1\left(\xi_{\s(3)}\right)}_\g,
\end{eqnarray*}
where we have used $(\E_2)$ (i.e., $\p_\g\circ f_2 -f_1\circ\ell_2=T_2$) in the second step.
As $\p_\g$ is a derivation for $\LB_\g$ and using the fact that $\p_\g\circ f_1=0$, one obtains:
$$
\lb{\p_\g\left(f_2\left(\xi_{\s(1)},\xi_{\s(2)}\right)\right)\!,f_1\left(\xi_{\s(3)}\right)}_\g = \p_\g\left(\lb{f_2\left(\xi_{\s(1)},\xi_{\s(2)}\right),f_1\left(\xi_{\s(3)}\right)}_\g\right).
$$ 
Finally, because $p\circ \p_\g = 0$,
\begin{eqnarray*}
\lefteqn{\ds -\rna{0.7}\sum_{\s\in S_{2,1}} \chi(\s)\, p\circ T_2\left(\ell_2\left(\xi_{\s(1)},\xi_{\s(2)}\right)\!,\xi_{\s(3)}\right)=}\\
&& \ds p\left(\rna{0.7}\sum_{\s\in S_{2,1}} \chi(\s) \lb{\lb{f_1\left(\xi_{\s(1)}\right),f_1\left(\xi_{\s(2)}\right)}_\g,f_1\left(\xi_{\s(3)}\right)}_\g\right)=0,
\end{eqnarray*}
where we have used the graded Jacobi identity satisfied by $\LB_\g$, to obtain the last line. This shows that the map $\ell_2$ satisfies ($\J_3$).
\begin{rem}\label{rem:ell2}
The skew-symmetric graded linear map $\ell_2$ of degree $0$ which satisfies ($\E_2$) is unique and given by (\ref{eq:ell2}). Using (\ref{eq:ell2}) and the definition of $T_2$, we obtain that, for all $\xi_1,\xi_2\in H(\g,\p_\g)$,
$$
\ell_2(\xi_1,\xi_2) =-p\circ T_2(\xi_1,\xi_2) = p\left(\lb{f_1(\xi_1),f_1(\xi_2)}_\g \right).
$$
In other words, the map $\ell_2 : \bigwedge\nolimits^2 H(\g,\p_\g)\to H(\g,\p_\g)$ is the map induced by the graded Lie bracket $\LB_\g$ on $H(\g,\p_\g)$. For this reason, we sometimes denote $\ell_2$ also by $\LB_\g$.
\end{rem}
%

 Let us now prove the part (b) of the proposition. To do this, we suppose that $m\geq 3$ and that  $f_2,\dots,f_{m-1}$ and $\ell_2,\dots,\ell_{m-1}$ are skew-symmetric graded linear maps, of degrees $\deg(\ell_k)=2-k$ and $\deg(f_k)=1-k$, which satisfy the equations $(\J_2)$ -- $(\J_m)$ and $(\E_2)$ -- $(\E_{m-1})$. Then, we show (\textit{Step 1}), that 
$$
\p_\g\left(T_m(\xi_1,\dots,\xi_m)\right) = 0, \hbox{ for all } \xi_1,\dots,\xi_m\in H(\g,\p_\g).
$$
This indeed implies, by (\ref{eq:decomp-f1}), that the cocycle $T_m(\xi_1,\dots,\xi_m)$ decomposes as a coboundary (element in the image of $\p_\g$) plus an element in the image of $f_1$, which leads to the existence of both maps $f_m$ and $\ell_m$, satisfying the equation ($\E_{m}$). 

Then (\textit{Step 2}), we show that the obtained map $\ell_m$, satisfying ($\E_{m}$), necessarily also satisfies the equation ($\J_{m}$). 

\bigskip

(b) - \textit{ \underline{Step 1.}} Let $\xi_1,\dots,\xi_m\in H(\g,\p_\g)$ be homogeneous elements. Recall that we have:
\begin{equation}\label{eq:TSU}
T_m(\xi_1,\dots,\xi_m) = S_m(\xi_1,\dots,\xi_m) - U_m(\xi_1,\dots,\xi_m),
\end{equation}
where we define, for all $n\in\N^*$, and all $\zeta_1,\dots,\zeta_n\in H(\g,\p_\g)$:
\begin{equation}\label{S(m)}
   \renewcommand{\arraystretch}{1.7}
     \begin{array}{l}
\quad S_n(\zeta_1,\dots,\zeta_n) :=\\
\!\!\!
\ds\sum_{\renewcommand{\arraystretch}{0.5}
     \begin{array}{c}
          \scriptstyle j+k=n+1\\
          \scriptstyle j,k\geq2
     \end{array}}\!\!
 \sum_{\s\in S_{k,n-k}}\!\! \chi(\s)\,(-1)^{k(j-1)}f_j\left(\ell_k\left(\zeta_{\s(1)},\dots,\zeta_{\s(k)}\right)\!,\zeta_{\s(k+1)},\dots,\zeta_{\s(n)}\right)
     \end{array}
\end{equation}
 and
\begin{equation}\label{U(m)}
   \renewcommand{\arraystretch}{1.7}
     \begin{array}{l}
\quad U_n(\zeta_1,\dots,\zeta_n) :=\\
\!\!\!
\ds\sum_{\renewcommand{\arraystretch}{0.5}
    \begin{array}{c}
        \scriptstyle s+t=n\\
        \scriptstyle s,t\geq1\\
     \end{array}}
\!\!\!\!\sum_{\renewcommand{\arraystretch}{0.5}
    \begin{array}{c}
        \scriptstyle \tau\in S_{s,n-s}\\
        \scriptstyle \tau(1)=1
     \end{array}} \!\!\! \!\!\! \!
  \chi(\tau) \,e_{s,t}(\tau)
      \lb{f_s \left(\zeta_{\tau(1)},\dots,\zeta_{\tau(s)}\right), f_t \left(\zeta_{\tau(s+1)},\dots,\zeta_{\tau(n)}\right)}_\g,
     \end{array}
\end{equation}
with $e_{s,t}(\tau) =(-1)^{s-1}\cdot(-1)^{(t-1)\left(\sum\limits_{p=1}^s |\zeta_{\tau(p)}|\right)}$ and where $\chi(\s)$ (respectively, $\chi(\tau)$) stands for $\chi(\s;\zeta_1,\dots,\zeta_n)$ (respectively, $\chi(\tau;\zeta_1,\dots,\zeta_n)$).
For $j=2,\dots, m-1$, the equation ($\E_j$) can be written as $\p_\g\circ f_j = T_j + f_1\circ \ell_j$, so that
\begin{equation*}
   \renewcommand{\arraystretch}{1.7}
     \begin{array}{l}
\quad \p_\g\left(S_m(\xi_1,\dots,\xi_m)\right)=\\
\!\!\!
\ds\sum_{\renewcommand{\arraystretch}{0.5}
     \begin{array}{c}
          \scriptstyle j+k=m+1\\
          \scriptstyle j,k\geq2
     \end{array}}\!\!
 \sum_{\s\in S_{k,m-k}}\!\! \chi(\s)\,(-1)^{k(j-1)}T_j\left(\ell_k\left(\xi_{\s(1)},\dots,\xi_{\s(k)}\right)\!,\xi_{\s(k+1)},\dots,\xi_{\s(m)}\right)\\
  +\;\ds f_1\left(J_m(\xi_1,\dots,\xi_m)\right)=\\
 \!\!\!
\ds\sum_{\renewcommand{\arraystretch}{0.5}
     \begin{array}{c}
          \scriptstyle j+k=m+1\\
          \scriptstyle j,k\geq2
     \end{array}}\!\!
 \sum_{\s\in S_{k,m-k}}\!\! \chi(\s)\,(-1)^{k(j-1)}T_j\left(\ell_k\left(\xi_{\s(1)},\dots,\xi_{\s(k)}\right)\!,\xi_{\s(k+1)},\dots,\xi_{\s(m)}\right),
     \end{array}
\end{equation*}
 where we have used the equation ($\J_m$) (in the case $\ell_1=0$, see (\ref{eq:Jn})), in the second step.
 Now, using the writing of $T_j$, for $2\leq j\leq m-1$, we get:
 \begin{equation*}
  \p_\g\left(S_m(\xi_1,\dots,\xi_m)\right) = \al_m(\xi_1,\dots,\xi_m) + \bl_m(\xi_1,\dots,\xi_m) + \cl_m(\xi_1,\dots,\xi_m),
 \end{equation*}
 where, for all $n\in\N^*$ and all homogeneous $\zeta_1,\dots,\zeta_n\in H(\g,\p_\g)$, we have defined:
\begin{equation*}
   \renewcommand{\arraystretch}{1.7}
     \begin{array}{l}
\quad \al_n(\zeta_1,\dots,\zeta_n) :=\\
\!\!\!\ds
\sum_{\renewcommand{\arraystretch}{0.5}
     \begin{array}{c}
          \scriptstyle p+q+k=n+2\\
          \scriptstyle p,q,k\geq2
     \end{array}}\!\!\!
 \sum_{\renewcommand{\arraystretch}{0.9}
     \begin{array}{c}
          \scriptstyle \a\in S_{q-1,p-1}^{k+1}\\
          \scriptstyle \s\in S_{k,n-k}
     \end{array}} \!\! 
 \chi(\s;\zeta_1,\dots,\zeta_n)\, \chi(\a;\zeta_{\s(k+1)},\dots,\zeta_{\s(n)})\cdot(-1)^{k(p+q) + q(p-1)}\cdot\\
 f_p\left(\ell_q\left(\ell_k\left(\zeta_{\s(1)},\dots,\zeta_{\s(k)}\right),\zeta_{\s\a(k+1)},\dots,\zeta_{\s\a(k+q-1)}\right),\zeta_{\s\a(k+q)},\dots,\zeta_{\s\a(n)}\right),
 \end{array}
 \end{equation*}
 and
 \begin{equation*}
   \renewcommand{\arraystretch}{1.7}
 \begin{array}{l}
 \quad \bl_n(\zeta_1,\dots,\zeta_n) :=\\
\!\!\!\ds
\sum_{\renewcommand{\arraystretch}{0.5}
     \begin{array}{c}
          \scriptstyle p+q+k=n+2\\
          \scriptstyle p,q,k\geq2
     \end{array}}\!\!
  \sum_{\renewcommand{\arraystretch}{0.9}
     \begin{array}{c}
          \scriptstyle \a\in S_{q,p-2}^{k+1}\\
          \scriptstyle \s\in S_{k,n-k}
     \end{array}}\!\!
 \chi(\s;\zeta_1,\dots,\zeta_n) \,\chi(\a;\zeta_{\s(k+1)},\dots,\zeta_{\s(n)}) \cdot\\
 \,(-1)^{k(p+q)}(-1)^{q(p-1)} \cdot 
 (-1)^{q+\left(\sum\limits_{r=1}^{k}|\zeta_{\s(r)}| + k\right)\cdot\left(\sum\limits_{s=k+1}^{k+q}|\zeta_{\s\a(s)}|\right)}\cdot\\
 f_p\left(\ell_q\left(\zeta_{\s\a(k+1)},\dots,\zeta_{\s\a(k+q)}\right),\ell_k\left(\zeta_{\s(1)},\dots,\zeta_{\s(k)}\right),\zeta_{\s\a(k+q+1)},\dots,\zeta_{\s\a(n)}\right),
     \end{array}
\end{equation*}
 and finally
 \begin{equation*}
   \renewcommand{\arraystretch}{1.7}
 \begin{array}{l}
 \quad \cl_n(\zeta_1,\dots,\zeta_n) :=\\
-\!\!\!\ds
\sum_{\renewcommand{\arraystretch}{0.5}
     \begin{array}{c}
          \scriptstyle j+k=n+1\\
          \scriptstyle j,k\geq2
     \end{array}}\!\!
 \sum_{\s\in S_{k,n-k}}\!\! 
 \sum_{\renewcommand{\arraystretch}{0.5}
     \begin{array}{c}
          \scriptstyle a+b=j\\
          \scriptstyle a,b\geq 1
     \end{array}}\!\!
  \sum_{\b\in S_{a-1,b}^{k+1}}\!\! 
 \chi(\s;\zeta_1,\dots,\zeta_n) \,\chi(\b;\zeta_{\s(k+1)},\dots,\zeta_{\s(n)}) \cdot\\
 \, (-1)^{k(j-1)}(-1)^{a-1} \cdot
 (-1)^{(b-1)\left(\sum\limits_{r=1}^{k}|\zeta_{\s(r)}| + k + \sum\limits_{s=k+1}^{k+a-1}|\zeta_{\s\b(s)}|\right)}
 \cdot\\
\lb{f_a\left(\ell_k\left(\zeta_{\s(1)},\dots,\zeta_{\s(k)}\right),\zeta_{\s\b(k+1)},\dots,\zeta_{\s\b(k+a-1)}\right),f_b\left(\zeta_{\s\b(k+a)},\dots,\zeta_{\s\b(n)}\right)}_\g.
     \end{array}
\end{equation*}
 Here, for $r,s,t\in\N$, we have denoted by $S_{s,t}^{r+1}$ the set of all the permutations $\s$ of $\lbrace r+1, \dots, r+s+t\rbrace$, such that $\s(r+1)<\cdots<\s(r+s)$ and $\s(r+s+1)<\cdots<\s(r+s+t)$. 
 A permutation  $\s\in S_{s,t}^{r+1}$ can also be seen as a permutation of $\lbrace 1, \dots, r+s+t\rbrace$, simply by fixing 
 ${\s}_{\vert_{\lbrace 1, \dots, r\rbrace}} = {\id}_{\vert_{\lbrace 1, \dots, r\rbrace}}$. 
 
%
\begin{rem}\label{just1}
Let us justify how one obtains that the sum 
\begin{equation*}
   \renewcommand{\arraystretch}{1.7}
     \begin{array}{l}
     \dl(\xi_1,\dots,\xi_m):=\\
\ds\sum_{\renewcommand{\arraystretch}{0.5}
     \begin{array}{c}
          \scriptstyle j+k=m+1\\
          \scriptstyle j,k\geq2
     \end{array}}\!\!
 \sum_{\s\in S_{k,m-k}}\!\! \chi(\s)\,(-1)^{k(j-1)}T_j\left(\ell_k\left(\xi_{\s(1)},\dots,\xi_{\s(k)}\right)\!,\xi_{\s(k+1)},\dots,\xi_{\s(m)}\right)
 \end{array}
 \end{equation*}
 is given by 
 $
 \al_m(\xi_1,\dots,\xi_m) + \bl_m(\xi_1,\dots,\xi_m) + \cl_m(\xi_1,\dots,\xi_m),
 $
using only the definition of the $T_j$.
Let $\xi_1,\dots,\xi_m\in H(\g,\p_\g)$ be homogeneous elements and let $j,k\geq 2$ with $ j+k=m+1$, and $\s\in S_{k,m-k}$. In order to simplify the notation, we denote by $\eta_1:=\ell_k\left(\xi_{\s(1)},\dots,\xi_{\s(k)}\right)$ and $\eta_2 := \xi_{\s(k+1)},\dots,\eta_j := \xi_{\s(m)}$ and write:
\begin{equation*}
   \renewcommand{\arraystretch}{1.7}
     \begin{array}{l}
T_j\left(\eta_1,\eta_2,\dots,\eta_j\right)=\\
\!\!\!\!\!\!\!\ds
 \sum_{\renewcommand{\arraystretch}{0.5}
     \begin{array}{c}
          \scriptstyle p+q=j+1\\
          \scriptstyle p,q\geq2
     \end{array}}\!\!
  \sum_{\gamma\in S_{q,j-q}}\!\! \!
\chi(\gamma;\eta_{1},\dots,\eta_{j})  (-1)^{q(p-1)} f_p\left(\ell_q\!\left(\eta_{\gamma(1)},\dots,\eta_{\gamma(q)}\right)\!\!,\eta_{\gamma(q+1)},\dots,\eta_{\gamma(j)}\right)\\
 %
 %
 - \!\!\!\ds
 \sum_{\renewcommand{\arraystretch}{0.5}
     \begin{array}{c}
          \scriptstyle a+b=j\\
          \scriptstyle a,b\geq 1
     \end{array}}\!\!
  \sum_{\renewcommand{\arraystretch}{0.5}
     \begin{array}{c}
          \scriptstyle \gamma'\in S_{a,j-a}\\
          \scriptstyle \gamma'(1) =1
     \end{array}}\!\!\!\!
 \chi(\gamma';\xi_{\s(k+1)},\dots,\xi_{\s(m)}) 
 (-1)^{a-1+(b-1)\left(\sum\limits_{r=1}^{a}|\eta_{\gamma'(r)}| \right)}\cdot\\
\qquad \qquad \qquad \qquad \qquad \qquad \qquad \lb{f_a\left(\eta_{\gamma'(1)},\dots,\eta_{\gamma'(a)}\right),f_b\left(\eta_{\gamma'(a+1)},\dots,\eta_{\gamma'(j)}\right)}_\g.
 \end{array}
 \end{equation*}
Then, the second sum leads easily to $\cl_m(\xi_1,\dots,\xi_m)$ and for the first sum, one has to separate the two cases where the permutation $\gamma\in S_{q,j-q}$, which appears in the sum, satisfies $\gamma(1)=1$ or $\gamma(q+1)=1$, to obtain respectively the terms $\al_m(\xi_1,\dots,\xi_m)$ and $\bl_m(\xi_1,\dots,\xi_m)$. Indeed, if $\gamma(1)=1$, then there exists $\a\in S_{q-1,p-1}^{k+1}$ such that: 
$$
\begin{array}{rcl}
\eta_{\gamma(1)} &=& \ell_k(\xi_{\s(1)}, \dots, \xi_{\s(k)}),\\
\eta_{\gamma(2)} &=& \xi_{\s\a(k+1)},\\
&\vdots&\\
\eta_{\gamma(q)} &=& \xi_{\s\a(k+q-1)},
\end{array}
\hfill \qquad
\begin{array}{rcl}
&&\\
\eta_{\gamma(q+1)} &=& \xi_{\s\a(k+q)}\\
&\vdots&\\
\eta_{\gamma(j)} &=&  \xi_{\s\a(m)}.
\end{array}
$$
By checking that $\chi(\gamma;\eta_1,\dots,\eta_j) = \chi(\a;\xi_{\s(k+1)},\dots,\xi_{\s(m)})$, one obtains the sum $\al_m(\xi_1,\dots,\xi_m)$.
 In the case $\gamma(q+1)=1$, one can rather write:
$$
\begin{array}{rcl}
\eta_{\gamma(1)} &=& \xi_{\s\a(k+1)},\\
&\vdots&\\
\eta_{\gamma(q)} &=& \xi_{\s\a(k+q)},\\
\eta_{\gamma(q+1)} &=& \ell_k(\xi_{\s(1)}, \dots, \xi_{\s(k)}),
\end{array}
\hfill \qquad
\begin{array}{rcl}
\eta_{\gamma(q+2)} &=& \xi_{\s\a(k+q+1)},\\
&\vdots&\\
\eta_{\gamma(j)} &=&  \xi_{\s\a(m)},\\
&&
\end{array}
$$
with $\a\in S_{q,p-2}^{k+1}$. It is then possible to compute that $\sgn(\gamma)=\sgn(\a)\cdot (-1)^{q}$ and 
$\e(\gamma;\eta_1,\dots,\eta_j) = \e(\a;\xi_{\s(k+1)},\dots,\xi_{\s(m)})\cdot (-1)^{\left(\sum\limits_{s=1}^k|\xi_{\s(s)}| + k\right)\cdot\left(\sum\limits_{r=k+1}^{k+q}|\xi_{\s\a(r)}|\right)}$. This permits one to obtain the sum $\bl_m(\xi_1,\dots,\xi_m)$.
 \end{rem}
%

 Now, we will successively show that both sums $\al_m(\xi_1,\dots,\xi_m)$ and $\bl_m(\xi_1,\dots,\xi_m)$ are equal to zero. 
 To do this, we prove the following lemmas.
 \begin{lma}\label{lma:al_n}
Let $n\in\N^*$. Suppose that the equations ($\J_j$) for $1\leq j\leq n-1$ are satisfied by the maps $\ell_1=0,\ell_2, \dots,\ell_{n-1}$, then 
$$
\al_n(\zeta_1,\dots,\zeta_n) =0, \quad \hbox{ for all }\; \zeta_1,\dots,\zeta_n\in H(\g,\p_\g).
$$
 \end{lma}
 \begin{proof}[proof of lemma \ref{lma:al_n}]
 Let $\zeta_1,\dots,\zeta_n\in H(\g,\p_\g)$.
 For $p,q,k\geq2$ such that $p+q+k=n+2$, and for $\s\in S_{k, n-k}$ and $\a\in S_{q-1,p-1}^{k+1}$, the permutation $\s\circ\a\in\S_{n}$ can be uniquely written as $\s\circ\a = \rho\circ\beta$, with $\rho\in S_{n-p+1,p-1}$ and $\b\in S_{k,q-1}$. Using this, one obtains:
\begin{equation*}
   \renewcommand{\arraystretch}{1.7}
     \begin{array}{l}
\al_n(\zeta_1,\dots,\zeta_n) =\\
\ds
\sum_{p=2}^{n-2}
 \sum_{\rho\in S_{n-p+1,p-1}}\!\!\!\!\!
 \chi(\rho) (-1)^{(n-p)(p-1)}
f_p\left(J_p(\zeta_{\rho(1)},\dots,\zeta_{\rho(n-p+1)}), \zeta_{\rho(n-p+2)},\dots,\zeta_{\rho(n)}\right)\!,
 \end{array}
 \end{equation*}
 where $\chi(\rho)$ stands for $\chi(\rho;\zeta_1,\dots,\zeta_n)$ and  $J_p$ is defined in (\ref{eq:Jn}). For every $2\leq p\leq n-2$ and every $\rho\in S_{n-p+1,p-1}$, one has $J_p(\zeta_{\rho(1)},\dots,\zeta_{\rho(n-p+1)})=0$, by ($\J_{n-p+1}$), where $n-p+1 = k+q-1$ runs through all integers between $3$ and $n-1$. 
Hence $\al_n(\zeta_1,\dots,\zeta_n)=0$.
 \end{proof}
According to this lemma, and because the maps $\ell_1=0,\ell_2,\dots,\ell_{m-1}$ are supposed to satisfy the equations ($\J_1$) -- ($\J_{m-1}$), we have $\al_m(\xi_1,\dots,\xi_m)=0$.
Let us now consider the sum $\bl_m(\xi_1,\dots,\xi_m)$. It is also zero, according to the following:
\begin{lma}\label{lma:bl_n}
Let $n\in\N^*$.  For all $\zeta_1,\dots,\zeta_n\in H(\g,\p_\g)$, we have
$$
\bl_n(\zeta_1,\dots,\zeta_n) =0.
$$
\end{lma}
\begin{proof}[proof of lemma \ref{lma:bl_n}]
This result follows from the skew-symmetry of the maps $f_1,\dots, f_n$, making the sum $\bl_n(\zeta_1,\dots,\zeta_n)$ equal to minus itself.
 \end{proof}

 Now, we consider the term $\p_\g\left(U_m(\xi_1,\dots,\xi_m)\right)$.
 As $\p_\g$ is a graded derivation for $\LB_\g$ and because $\LB_\g$ is skew-symmetric, one has, for all $\zeta_1,\dots,\zeta_m\in H(\g,\p_\g)$ and all $s,t\in \lbrace 1,\dots m-1\rbrace$ such that $s+t=m$:
 \begin{eqnarray*}
\lefteqn{\p_\g\left(\lb{f_s\left(\zeta_1,\dots,\zeta_s\right),f_t\left(\zeta_{s+1},\dots,\zeta_m\right)}_\g\right) =}\\
&& \lb{\p_\g\left(f_s\left(\zeta_1,\dots,\zeta_s\right)\right),f_t\left(\zeta_{s+1},\dots,\zeta_m\right)}_\g \\
&&- (-1)^{|f_s\left(\zeta_1,\dots,\zeta_s\right)|\left(1 +  |\p_\g\left(f_t\left(\zeta_{s+1},\dots,\zeta_m\right)\right)|\right)} 
\lb{\p_\g\left(f_t\left(\zeta_{s+1},\dots,\zeta_m\right)\right), f_s\left(\zeta_1,\dots,\zeta_s\right)}_\g.
\end{eqnarray*}
 Using this, the one-to-one correspondence between the set $\lbrace \tau\in S_{s,m-s}\mid \tau(1)=1\rbrace$ and the set $\lbrace \tau'\in S_{t,m-t}\mid \tau'(t+1)=1\rbrace$ and finally the fact that $S_{s,m-s} =\lbrace \tau\in S_{s,m-s}\mid \tau(1)=1\rbrace \sqcup \lbrace \tau\in S_{s,m-s}\mid \tau(s+1)=1\rbrace $, we obtain that: 
 $$
    \renewcommand{\arraystretch}{1.7}
     \begin{array}{l}
\quad \p_\g\left(U_m(\xi_1,\dots,\xi_m)\right) =\\
\!\!\!
\ds\sum_{\renewcommand{\arraystretch}{0.5}
    \begin{array}{c}
        \scriptstyle s+t=m\\
        \scriptstyle s,t\geq1\\
     \end{array}}
\!\!\sum_{\tau\in S_{s,m-s}} \!\!\! 
  \chi(\tau) \, e_{s,t}(\tau)
      \lb{\p_\g\left(f_s\left(\xi_{\tau(1)},\dots,\xi_{\tau(s)}\right)\right),f_t\left(\xi_{\tau(s+1)},\dots,\xi_{\tau(m)}\right)}_\g.
     \end{array}
$$
Finally, it remains for $\p_\g\left(T_m(\xi_1,\dots,\xi_m)\right)$:
\begin{eqnarray*}
\lefteqn{ \p_\g\left(T_m(\xi_1,\dots,\xi_m)\right) = \p_\g\left(S_m(\xi_1,\dots,\xi_m)\right) - \p_\g\left(U_m(\xi_1,\dots,\xi_m)\right)}\\
&= & c_m(\xi_1,\dots,\xi_m)\\
&-&\!\!\!
\ds\sum_{\renewcommand{\arraystretch}{0.5}
    \begin{array}{c}
        \scriptstyle s+t=m\\
        \scriptstyle s,t\geq1\\
     \end{array}}
\!\!\sum_{\tau\in S_{s,m-s}} \!\!\! 
  \chi(\tau) \,e_{s,t}(\tau)
  \lb{\p_\g\left(f_s\left(\xi_{\tau(1)},\dots,\xi_{\tau(s)}\right)\right),f_t\left(\xi_{\tau(s+1)},\dots,\xi_{\tau(m)}\right)}_\g.
\end{eqnarray*}
We now point out that, for all $n\in\N^*$ and for all $\zeta_1,\dots, \zeta_n\in H(\g,\p_\g)$,
\begin{eqn}{eq:cn}
\lefteqn{c_n(\zeta_1,\dots, \zeta_n) =}\\
&&\!\!\!\!\!\!\!\!\! \!\!\!\!\!\!\!\!\!\ds\sum_{\renewcommand{\arraystretch}{0.5}
    \begin{array}{c}
        \scriptstyle s+t=n\\
        \scriptstyle s,t\geq1\\
     \end{array}}
\!\!\sum_{\tau\in S_{s,n-s}} \!\!\! 
  \chi(\tau) \,e_{s,t}(\tau)
\lb{(f_1\circ\ell_s + S_s)\left(\zeta_{\tau(1)},\dots,\zeta_{\tau(s)}\right),f_t\left(\zeta_{\tau(s+1)},\dots,\zeta_{\tau(n)}\right)}_\g.
\end{eqn}
We use once more the equation ($\E_s$) and (\ref{eq:TSU}) to write  $\p_\g\circ f_s = f_1\circ \ell_s + T_s = f_1\circ \ell_s + S_s - U_s$, for $s=1,\dots,m-1$, and to obtain:
\begin{eqnarray*}
\lefteqn{ \p_\g\left(T_m(\xi_1,\dots,\xi_m)\right) =}\\
&&\!\!\!\!\!\!\!\!\! \!\!\!\!\!\!\!\!\!\ds\sum_{\renewcommand{\arraystretch}{0.5}
    \begin{array}{c}
        \scriptstyle s+t=m\\
        \scriptstyle s,t\geq1\\
     \end{array}}
\!\!\sum_{\tau\in S_{s,m-s}} \!\!\! 
  \chi(\tau) \,e_{s,t}(\tau)
   \lb{U_s\left(\xi_{\tau(1)},\dots,\xi_{\tau(s)}\right),f_t\left(\xi_{\tau(s+1)},\dots,\xi_{\tau(m)}\right)}_\g.
\end{eqnarray*}
Written differently, this reads as follows:
\begin{equation}\label{eq:TRm}
\p_\g\left(T_m(\xi_1,\dots,\xi_m)\right) = R_m(\xi_1,\dots,\xi_m),
\end{equation}
where we have introduced the following notation (because we will need this notation later): 
for all $n\in\N^*$ and all $\zeta_1,\dots,\zeta_n\in H(\g,\p_\g)$,
\begin{eqnarray*}
\lefteqn{ R_n(\zeta_1,\dots,\zeta_n) :=}\\
&&\!\!\!\!\!\!\!\!\! \!\!\!\!\!\!\!\!\!\ds\sum_{\renewcommand{\arraystretch}{0.5}
    \begin{array}{c}
        \scriptstyle a+b+t=n\\
        \scriptstyle a,b,t\geq1\\
     \end{array}}
\!\!\sum_{\renewcommand{\arraystretch}{0.7}
    \begin{array}{c}
        \scriptstyle \tau\in S_{a+b,t}\\
        \scriptstyle \s\in S_{a,b}\\
        \scriptstyle \s(1)=1\\
     \end{array}} \!\!\! 
  \chi(\tau;\zeta_1,\dots,\zeta_n)\;\chi(\s;\zeta_{\tau(1)},\dots,\zeta_{\tau(a+b)}) \;e_{a+b,t}(\tau) \; e_{a,b}(\tau\circ\s)\cdot\\
 &&\!\!\!\!\!\!\!\!\! \!\!\!\!\!\!\!\!\!     \lb{\lb{f_a\left(\zeta_{\tau\s(1)},\dots,\zeta_{\tau\s(a)}\right),f_b\left(\zeta_{\tau\s(a+1)},\dots,\zeta_{\tau\s(a+b)}\right)}_\g,
      f_t\left(\zeta_{\tau(a+b+1)},\dots,\zeta_{\tau(n)}\right)}_\g.
\end{eqnarray*}
It is then possible to show that this is zero, using the graded Jacobi identity satisfied by $\LB_\g$. Because we will need this result in another context, we show the following:
\begin{lma}\label{lma:R_m}
For $n\in\N^*$ and all $\zeta_1,\dots,\zeta_n\in H(\g,\p_\g)$, one has:
$$
R_n(\zeta_1,\dots,\zeta_n) =0.
$$
\end{lma}
\begin{proof}[proof of lemma \ref{lma:R_m}]
Let $\zeta_1,\dots,\zeta_n\in H(\g,\p_\g)$.
One first can show that 
\begin{eqnarray*}
\lefteqn{ 2\,R_n(\zeta_1,\dots,\zeta_n)=}\\
 &&\!\!\!\!\!\!\!\!\! \!\!\!\!\!\!\!\!\!  \ds\sum_{\renewcommand{\arraystretch}{0.5}
    \begin{array}{c}
        \scriptstyle a+b+t=n\\
        \scriptstyle a,b,t\geq1\\
     \end{array}}
\!\!\sum_{\rho\in S_{a,b,t}}
  \chi(\rho;\zeta_1,\dots,\zeta_n)\;
 e_{a+b,t}(\rho) \; e_{a,b}(\rho) \cdot\\
 &&\!\!\!\!\!\!\!\!\! \!\!\!\!\!\!\!\!\!        \lb{\lb{f_a\left(\zeta_{\rho(1)},\dots,\zeta_{\rho(a)}\right),f_b\left(\zeta_{\rho(a+1)},\dots,\zeta_{\rho(a+b)}\right)}_\g,
      f_t\left(\zeta_{\rho(a+b+1)},\dots,\zeta_{\rho(n)}\right)}_\g,
\end{eqnarray*}
where for $a,b,t\in\N$, $S_{a,b,t}$ is the set of all the permutations $\s\in\S_{a+b+t}$ of $\lbrace 1,\dots,a+b+t\rbrace$, satisfying:
$\s(1)<\cdots< \s(a)$, $\s(a+1)<\cdots< \s(a+b)$ and $\s(a+b+1)<\cdots< \s(a+b+t)$.
It is now possible to check that one has:
\begin{eqnarray*}
\lefteqn{ 6\,R_n(\zeta_1,\dots,\zeta_n) =
\ds\sum_{\renewcommand{\arraystretch}{0.5}
    \begin{array}{c}
        \scriptstyle a+b+t=n\\
        \scriptstyle a,b,t\geq1\\
     \end{array}}
\!\!\sum_{\rho\in S_{a,b,t}}
  \chi(\rho)\,(-1)^{e}\cdot}\\
&&\!\!\!\!\!\!\!\!\! \!\!\!\!\!\!\!\!\!   \Jac_\g\left(f_a\left(\zeta_{\rho(1)},\dots,\zeta_{\rho(a)}\right),f_b\left(\zeta_{\rho(a+1)},\dots,\zeta_{\rho(a+b)}\right),f_t\left(\zeta_{\rho(a+b+1)},\dots,\zeta_{\rho(n)}\right)\right),
\end{eqnarray*}
where $e\in\Z$ is an integer depending on $\zeta_1,\dots,\zeta_n$ and on the permutation $\rho$, and where, for all $x,y,z\in \g$, 
$$
\Jac_\g(x,y,z) := (-1)^{|x||z|}\lb{\lb{x,y}_g,z}_\g + (-1)^{|y||x|}\lb{\lb{y,z}_g,x}_\g + (-1)^{|z||y|}\lb{\lb{z,x}_g,y}_\g,
$$
 which is zero because of the graded Jacobi identity satisfied by $\LB_\g$. We now conclude that
$R_n(\zeta_1,\dots,\zeta_n) = 0$.
\end{proof}
This lemma, together with (\ref{eq:TRm}), imply that $\p_\g\left(T_m(\xi_1,\dots,\xi_m)\right) = 0$.
This fact means that, for all $\xi_1,\dots,\xi_m\in H(\g,\p_\g)$, the element $T_m(\xi_1,\dots,\xi_m)$ is a cocycle for the cochain complex $(\g,\p_\g)$. This allows us to define a skew-symmetric graded linear map $\ell_m: \bigwedge\nolimits^m H(\g,\p_\g) \to H(\g,\p_\g)$, of degree $2-m$, with the following formula:
\begin{equation}\label{eq:ell_m}
\ell_m(\xi_1,\dots,\xi_m) := -p\circ T_m(\xi_1,\dots,\xi_m),
\end{equation}
for all $\xi_1,\dots,\xi_m\in H(\g,\p_\g)$. As in the case $m=2$ and according to (\ref{eq:property-f1}), we also have the existence of a skew-symmetric graded linear map $f_m: \bigwedge\nolimits^m H(\g,\p_\g) \to \g$, of degree $1-m$, which satisfies the equation ($\E_m$):
$$
T_m(\xi_1,\dots,\xi_m)= \p_\g\left(f_m(\xi_1,\dots,\xi_m)\right) - f_1\left(\ell_m(\xi_1,\dots,\xi_m)\right),
$$
for all $\xi_1,\dots,\xi_m\in H(\g,\p_\g)$.

\bigskip

(b) - \textit{ \underline{Step 2.}} It remains to show, using the equations ($\J_1$) -- ($\J_m$) and ($\E_1$) -- ($\E_{m-1}$), satisfied by the maps $\ell_1,\dots,\ell_{m-1}$ and $f_1,\dots,f_{m-1}$ and the equation ($\E_m$) also satisfied by the maps $\ell_m$ and $f_m$, that the map $\ell_m$, defined in (\ref{eq:ell_m}), satisfies necessarily, for all $\xi_1,\dots,\xi_{m+1}\in H(\g,\p_\g)$, the equation:
\begin{equation}\label{L(m+1)}\tag{$\J_{m+1}$}
\renewcommand{\arraystretch}{0.7} 
\sum_{\begin{array}{c}
  \scriptstyle j+k=m+2\\
  \scriptstyle j,k\geq 2
\end{array}}
\!\!\sum_{\s\in S_{k,m+1-k}} 
\chi(\s)(-1)^{k(j-1)}\, \ell_j\left(\ell_k\left(\xi_{\s(1)},\dots,\xi_{\s(k)}\right)\!,\xi_{\s(k+1)},\dots,\xi_{\s(m)}\right)=0.
\end{equation}
Let us fix $\xi_1,\dots,\xi_{m+1}\in H(\g,\p_\g)$.
By equations ($\E_1$) -- ($\E_m$), we know that the maps $\ell_j$, for $1\leq j\leq m$, can be written as $\ell_j = - p\circ T_j$.
Using the notation of the remark \ref{just1}, this implies that (\ref{L(m+1)}) is equivalent to:
\begin{equation*}
p\left(\dl_{m+1}(\xi_1,\dots,\xi_{m+1})\right)=0.
\end{equation*}
We also use the same reasoning as the one explained in the remark \ref{just1} to obtain:
\begin{equation*}
\dl_{m+1}(\xi_1,\dots,\xi_{m+1})=\left(\al_{m+1} + \bl_{m+1}+ \cl_{m+1}\right)(\xi_1,\dots,\xi_{m+1}).
\end{equation*}
Then, the lemma \ref{lma:al_n}, together with the fact that the maps $\ell_2,\dots,\ell_m$ satisfy the equations ($\J_j$) for $1\leq j\leq m$, imply that $\al_{m+1}(\xi_1,\dots,\xi_{m+1})=0$. Secondly, the lemma \ref{lma:bl_n} also says that $\bl_{m+1}(\xi_1,\dots,\xi_{m+1})=0$.
Finally it remains that:  
$$
\hbox{(\ref{L(m+1)}) is equivalent to: }\; p\left(\cl_{m+1}(\xi_1,\dots,\xi_{m+1})\right)=0,
$$
 which is also equivalent to say that  $\cl_{m+1}(\xi_1,\dots,\xi_{m+1})$ is a coboundary for the cochain complex $(\g,\p_\g)$. 
 As (\ref{eq:cn}) can be obtained without using anything but the definitions of $\cl_m$ and $S_s$, we also have:
\begin{eqnarray*}
\lefteqn{ \cl_{m+1}(\xi_1,\dots,\xi_{m+1})=}\\
&&\!\!\!\!\!\!\!\!\! \!\!\!\!\!\!\! \ds
\sum_{\renewcommand{\arraystretch}{0.5}
     \begin{array}{c}
          \scriptstyle p+q=m+1\\
          \scriptstyle q\geq 1, p\geq2
     \end{array}}\!\!\!
\sum_{\a\in S_{p,q}} \!\! \!\chi(\a) 
e_{p,q}(\a)\!\!\lb{\left(S_p\!+ \! f_1\circ \ell_p\right)\!\left(\xi_{\a(1)},\dots,\xi_{\a(p)}\right)\!\!,
 f_q\!\left(\xi_{\a(p+1)},\dots,\xi_{\a(m+1)}\right)}_\g\!.
\end{eqnarray*}
Now, we use $S_p=T_p + U_p$ and the equations ($\E_p$), satisfied by the maps $\ell_p$ and $f_p$, for $1\leq p\leq m$, to write
$S_p+ f_1\circ \ell_p=\p_\g\circ f_p + U_p$ and:
\begin{eqnarray*}
\lefteqn{ \cl_{m+1}(\xi_1,\dots,\xi_{m+1})=}\\
&&\!\!\!\!\!\!\!\!\! \!\!\!\!\!\!\!\!\! \ds
\sum_{\renewcommand{\arraystretch}{0.5}
     \begin{array}{c}
          \scriptstyle p+q=m+1\\
          \scriptstyle q\geq 1, p\geq2
     \end{array}}\!\!
\sum_{\a\in S_{p,q}}  \chi(\a) \;
 e_{p,q}(\a) \lb{\p_\g\left(f_p\left(\xi_{\a(1)},\dots,\xi_{\a(p)}\right)\right),
 f_q\left(\xi_{\a(p+1)},\dots,\xi_{\a(m+1)}\right)}_\g\\
&+& R_{m+1}(\xi_1,\dots,\xi_{m+1}).
\end{eqnarray*}
By lemma \ref{lma:R_m}, $R_{m+1}(\xi_1,\dots,\xi_{m+1})=0$, and using the bijection between $S_{p,q}$ and $S_{q,p}$, given by:
$$
\begin{array}{rcl}
S_{p,q} &\to& S_{q,p}\\
\a &\mapsto & \a' :=
   \left(
   \begin{array}{cccccc}
   \sc 1& \sc\cdots & \sc q & \sc q+1 & \sc\cdots & \sc p+q\\
   \sc\a(p+1) & \sc \cdots & \sc\a(p+q)& \sc\a(1)& \sc \cdots & \sc\a(p)
   \end{array}
   \right),
\end{array}
$$
for which 
\begin{eqnarray*}
\sgn(\a') &=& \sgn(\a)\cdot (-1)^{pq},\\
\e(\a';\xi_1,\dots,\xi_{p+q}) &=& \e(\a;\xi_1,\dots,\xi_{p+q})
\cdot (-1)^{\left(\sum\limits_{r=1}^{q}|\xi_{\a(r)}|\right)\cdot \left(\sum\limits_{r=q+1}^{p+q}|\xi_{\a(r)}|\right)},
\end{eqnarray*}
 and also using the skew-symmetry of $\LB_\g$ and the fact that $\p_\g$ is a graded derivation for $\LB_\g$, we finally obtain:
\begin{eqnarray*}
\lefteqn{2\,\cl_{m+1}(\xi_1,\dots,\xi_{m+1})=}\\
&&\!\!\!\!\!\!\!\!\! \!\!\!\!\!\!\!\!\!\ds
\sum_{\renewcommand{\arraystretch}{0.5}
     \begin{array}{c}
          \scriptstyle p+q=m+1\\
          \scriptstyle q, p\geq2
     \end{array}}\!\!
\sum_{\a\in S_{p,q}}  \chi(\a) \;
e_{p,q}(\a)\; \p_\g\left(\lb{f_p\left(\xi_{\a(1)},\dots,\xi_{\a(p)}\right),
              f_q\left(\xi_{\a(p+1)},\dots,\xi_{\a(m+1)}\right)}_\g\right).
 \end{eqnarray*}
 We have then obtained that $\cl_{m+1}(\xi_1,\dots,\xi_{m+1})$ is a coboundary for the cochain complex $(\g,\p_\g)$, so that 
 $ \p_\g\left(\cl_{m+1}(\xi_1,\dots,\xi_{m+1})\right) =0$
and the equation ($\J_{m+1}$) is satisfied. This finishes the proof of the proposition \ref{prop:magic}.
\end{proof}
\section{Deformations of Poisson structures via $L_\infty$-algebras}
\label{sec:defosPoisson}

In this section, we consider a family of dg Lie algebras, constructed from a family of Poisson structures in dimension three. We will then use the proposition~\ref{prop:magic}, to obtain a classification of all formal deformations of these Poisson structures in the generic case, together with an explicit formula for the representative of each equivalence classes of these deformations.

\subsection{Poisson structures in dimension three and their cohomology}
\label{subsec:jac}

In the following, $\A$ denotes the polynomial algebra in three generators $\A := \F[x,y,z]$, where $\F$ is an arbitrary field of characteristic zero. To each polynomial $\v\in\A$, one associates a Poisson structure $\PB_\v$ defined by:
\begin{equation}\label{eq:PBv}
\PB_\v := \pp{\v}{x}\, \pp{}{y}\wedge \pp{}{z} + \pp{\v}{y}\, \pp{}{z}\wedge \pp{}{x} + \pp{\v}{z}\, \pp{}{x}\wedge \pp{}{y}.
\end{equation}
In this context, the Poisson cohomology of $(\A,\PB_\v)$ is denoted by $H(\A,\PB_\v)$.
We also denote by $(\g_\v,\p_{\v},\LB_{S})$, the dg Lie algebra associated to the Poisson algebra $(\A,\PB_\v)$, as explained in the paragraph \ref{sbsec:dglie-poisson}. Notice that $\g_\v^k\simeq \lbrace 0\rbrace$, for all $k\geq 3$.
With these notations, and those of the previous section, we have:
$H^n(\g_\v,\p_{\v})=H^{n+1}(\A,\PB_\v)$, for all $n\in\Z$ (in fact, $n\in\N\cup \lbrace -1\rbrace$). As previously, for every cocycle $P$ of the cochain complex $(\g_\v,\p_\v)$, $\bar P$ denotes its cohomology class in $H(\g_\v,\p_\v)$.
As we want to use the result of the previous section (proposition \ref{prop:magic}), we need to choose representatives $\left(\V^n_k\right)_{k}$ of an $\F$-basis of $H^n(\g_\v,\p_{\v})$, for $n\in\Z$.
To do this, we use the results of \cite{AP1}, in which the polynomial $\v$ is supposed to be weight-homogeneous and with an isolated singularity (at the origin).
Let us recall that a polynomial $\v\in\F[x,y,z]$ is said to be \emph{weight homogeneous}
of (weighted) degree $\w(\v)\in\N$, if there exists
(unique) positive integers
$\w_1,\w_2,\w_3\in\N^*$ (the \emph{weights} of the variables $x$, $y$ and
$z$), without any common divisor, such that:
\begin{equation}\label{euler_form}
\w_1 \,x\, \pp{\v}{x} + \w_2\, y\, \pp{\v}{y} + \w_3\, z\, \pp{\v}{z}= \w(\v) \v.
\end{equation}
This equation is called the \emph{Euler Formula} and can also be written as: 
$\vec{e}_\w[\v] = \w(\v) \v$,
 where $\vec{e}_\w$ is the so-called \emph{Euler derivation} (associated to the weights of the variables), defined by: 
$$
\vec{e}_\w := \w_1 \,x\, \pp{}{x} + \w_2\, y\, \pp{}{y} + \w_3\, z\, \pp{}{z}.
$$
Recall that a weight homogeneous polynomial
$\v\in\F[x,y,z]$ is said to admit an \emph{isolated singularity} (at
the origin) if the vector space
\begin{eqnarray}
\label{A_sing}
\A_{sing}(\v):=\F[x,y,z]/\ideal{\pp{\v}{x}, \pp{\v}{y}, \pp{\v}{z}}
\end{eqnarray}
is finite-dimensional. Its dimension is then denoted by $\mu$ and
called the \emph{Milnor number} associated to $\v$.  When $\F=\C$,
this amounts, geometrically, to saying that the surface
$\Fe_{\v}:\{\v=0\}$ has a singular point only at the origin.

From now on, the polynomial $\v$ will always be a weight homogeneous
polynomial with an isolated singularity. The corresponding weights of the three variables
($\w_1$, $\w_2$ and $\w_3$)
are then fixed and the weight homogeneity of any polynomial in $\A=\F[x,y,z]$
has now to be understood as associated to these weights. In the following, $|\w|$ denotes the sum of the weights of the three variables $x$, $y$ and $z$: $|\w| := \w_1 + \w_2 + \w_3$ and we fix $u_0:=1,
u_1,\dots,u_{\mu-1}\in\A$, a family composed of weight homogeneous
polynomials in $\A$ whose images in $\A_{sing}(\v)$ give a basis of this
$\F$-vector space (and $u_0=1$). (For example, one can choose the polynomials $u_0, \dots,u_{\mu-1}$ as being monomials of $\F[x,y,z]$).

\begin{prp}[\cite{AP1}]\label{prp:ap}
Let $\v\in \A$ be a weight-homogeneous polynomial with an isolated singularity. Let $(\g_\v,\p_{\v},\LB_{S})$ denote the dg Lie algebra associated to the Poisson algebra $(\A,\PB_\v)$, as explained in the paragraph \ref{sbsec:dglie-poisson}, and where $\PB_\v$ is defined in (\ref{eq:PBv}). Here we give explicit representatives for $\F$-bases of the Poisson cohomology spaces associated to $(\A,\PB_\v)$ or equivalently to $(\g_\v,\p_{\v})$.
\begin{enumerate}
\item An $\F$-basis of the first cohomology space $H^{-1}(\g_\v,\p_{\v})=H^{0}(\A,\PB_\v)$ is given by:
$$
\bb^{-1}_\v := \left(\overline{\v^i},\; i\in\N\right);
$$
\item An $\F$-basis of the space $H^0(\g_\v,\p_{\v})=H^{1}(\A,\PB_\v)$ is given by:
 $$
 \bb^0_\v := \left\lbrace
  \begin{array}{cc}
  \ds (0) \quad &\hbox{ if } \w(\v) \not= |\w|,\\
  \ds  \left(\overline{\v^i\, \vec{e}_\w}, \;i\in\N\right)\quad& \hbox{ if } \w(\v)= |\w|;
  \end{array}
  \right.
 $$
\item An $\F$-basis of the space $H^1(\g_\v,\p_{\v})=H^{2}(\A,\PB_\v)$ is given by:
$$
\bb^1_\v := \left( \overline{\v^i\, u_q \PB_{\v}},\; i\in\N, q\in\E_\v\right) \cup \left(\overline{\PB_{u_r}}, \;1\leq r\leq \mu-1\right),
$$
where 
$$
\E_\v:=\left\lbrace \begin{array}{cc}
   \lbrace 1,\dots,\mu-1\rbrace \quad & \hbox{ if } \w(\v) \not= |\w|,\\
   \lbrace 0,\dots,\mu-1\rbrace \quad & \hbox{ if } \w(\v) = |\w|,
   \end{array}
   \right.
$$ 
and where the skew-symmetric biderivation $\PB_{u_q}$ is naturally obtained by replacing $\v$ by $u_q$ in (\ref{eq:PBv});
\item An $\F$-basis of the space $H^2(\g_\v,\p_{\v})=H^{3}(\A,\PB_\v)$ is given by:
$$
\bb^2_\v := \left(\overline{\v^i\, u_s \D},\; i\in\N,\, 0\leq s\leq \mu-1\right),
$$
where $\D$ is the skew-symmetric triderivation of $\A$, defined by: 
$$
\D := \pp{}{x} \wedge \pp{}{y} \wedge \pp{}{z};
$$

\item For $k\geq 3$, 
$$
H^k(\g_\v,\p_{\v})=H^{k+1}(\A,\PB_\v)\simeq\lbrace 0\rbrace.
$$
\end{enumerate}
\end{prp}
\begin{rem}
More precisely, the basis of $H^{2}(\A,\PB_\v)$ given here is obtained by using the proposition 4.8 and the equality (27) of \cite{AP1}.
\end{rem}
\subsection{A suitable quasi-isomorphism between $H(\g_\v,\p_\v)$ and $\g_\v$}
Similarly to the definition (\ref{eqf_1}), we now have a linear graded map $f_1^\v$ of degree $0$, associated to the bases $\bb^{-1}_\v, \bb^{0}_\v, \bb^{1}_\v, \bb^{2}_\v$:
\begin{equation}\label{eq:f1v}
  \begin{array}{lccl}
   f_1^\v :& H^\ell(\g_\v,\p_\v)&\to&Z^\ell(\g_\v,\p_\v) \\
   &\\
    & \xi = \sum_k\lambda_k^\ell \overline{\V^\ell_k} &\mapsto& \sum_k\lambda_k^\ell \V^\ell_k,
  \end{array}
 \end{equation}
where $\xi =\sum_k\lambda_k^\ell \overline{ \V^\ell_k}$ is the unique decomposition of $\xi$ in the basis $\bb^{\ell}_\v$, $\ell=-1,0,1,2$, for which the elements $\left(\V^\ell_k\right)_k$ denote here the representatives, chosen in the previous proposition \ref{prp:ap}, of the basis $\bb^{\ell}_\v$.

Using the proposition \ref{prop:magic} and the bases $\bb^{-1}_\v, \bb^{0}_\v, \bb^{1}_\v, \bb^{2}_\v$ of the Poisson cohomology spaces associated to $(\g_\v,\p_\v)$, we construct an \l-algebra structure on $H(\g_\v,\p_\v)$: $\left(H(\g_\v,\p_\v),\ell_1=0,\ell_2=\LB_S,\ell_3,\dots\right)$, and a (weak) \l-morphism
$$
f_\bullet^\v=\left(f_n^\v :  \bigoplus\nolimits^n H(\g_\v,\p_\v) \to \g_\v \right)_{n\in\N^*},
$$
which extends $f_1^\v$, thus is a quasi-isomorphism.
We indeed prove the following:

\begin{thm}\label{prop:defosLinfini}
Let $\v\in\A=\F[x,y,z]$ be a weight-homogeneous polynomial, with an isolated singularity and let $\PB_\v$ be the associated Poisson bracket defined in~(\ref{eq:PBv}). 
Let $(\g_\v,\p_\v,\LB_S)$ be the dg Lie algebra associated to the Poisson cohomology complex of $(\A, \PB_\v)$, as explained in the paragraph~\ref{sbsec:dglie-poisson}. 
For simplicity, we denote by $H_\v$ the space $H(\g_\v,\p_\v)$, and for all $i\in\N^*$, $H_\v^i$ the space $H^i(\g_\v,\p_\v)$ (the $i$-th cohomology space associated to $(\g_\v,\p_\v)$).
We fix $f_1^\v$ as being the map defined in (\ref{eq:f1v}) and $\ell_1^\v:H(\g,\p_\g) \to H(\g,\p_\g)$ as being the trivial map.
We also fix the map $\ell_2^\v$ as being the bracket induced by the Schouten bracket $\LB_S$, i.e., 
$$
\ell_2^\v(\overline{x}, \overline{y}) :=\overline{ \lb{x,y}_S}, \quad \hbox{ for all }\; x,y\in\g_\v.
$$

There exist an \l-algebra structure on $H_\v :=H(\g_\v,\p_\v)$, denoted by $\ell_\bullet^\v :=  \left(\ell_i^\v \right)_{i\in\N^*}$ (with $\ell_1^\v$ and $\ell_2^\v$ given previously) and a quasi-isomorphism $f_\bullet^\v :=  \left(f_i^\v \right)_{i\in\N^*}$ (extending $f_1^\v$) from $H_\v$ to the dg Lie algebra $(\g_\v,\p_\v,\LB_S)$, satisfying the following properties:
\begin{enumerate}
\item[$(P_1)$] The map $f_2^\v$ is defined by the values given in the table \ref{tab:1}, for the case $\w(\v)\not=|\w|$, and in the table \ref{tab:2}, for the case $\w(\v)=|\w|$;
\item[$(P_2)$] For all $i\geq 2$, the map $\ell_i^\v$ is zero on $H^1_\v$:
$$
{\ell_i^\v}_{\vert_{\left(H_\v^1\right)^{\otimes i}}}=0,\; \hbox{ for all } i\geq 2;
$$ 
\item[$(P_3)$] For all $i\geq 3$, the map $f_i^\v$ is zero on $H^1_\v$:
$$
{f_i^\v}_{\vert_{\left(H_\v^1\right)^{\otimes i}}}=0,\; \hbox{ for all } i\geq 3.
$$
\end{enumerate}

\end{thm}
\begin{table}[ht]
  \caption{Case $\w(\v)\not=|\w|$. 
           The values of the linear map $f_2^\v$ on the elements of the bases $\bb^{i}_\v$ and $\bb^{j}_\v$ of the spaces $H_\v^i$ and $H_\v^j$, for $i,j=-1, 1, 2$. 
                 Notice that in this case, $H^0_\v = \lbrace 0\rbrace$. 
                 In this table, $F(\v), G(\v)$ are arbitrary elements of $\F[\v]$ and $1\leq k,l,s,t\leq \mu-1$.}
                 \label{tab:1}
  \begin{center}
  \renewcommand{\arraystretch}{1.5}
  \begin{tabular}{|c|c|c|}
    \hline 
    $H^i_\v \times H^j_\v$ & $(\overline{\V^i},\,\overline{\V^j})\in \bb^{i}_\v\times \bb^{j}_\v$ & $f_2^\v(\overline{\V^i},\overline{\V^j})\in \g_\v^{i+j-1}$\\
    \hline
    \hline
       $H^{-1}_\v \times H^{1}_\v$   & $\begin{array}{c}
                                                               \left(\overline{F(\v)}, \,\overline{G(\v)\, u_l \PB_\v}\right)\\
                                                                \left(\overline{F(\v)},\, \overline{\PB_{u_s}}\right)
                                                                  \end{array}$ & $\begin{array}{c}0\\F'(\v)\, u_s\end{array}$  \\
      \hline
      $H^{-1}_\v \times H^{2}_\v$    & 
                $ \begin{array}{c}
                    \left(\overline{F(\v)}, \,\overline{G(\v)\, u_l \D}\right)\\
                    \left(\overline{F(\v)}, \,\overline{G(\v)\, \D}\right)
                   \end{array}$ &   
                $\begin{array}{c}
                    0\\ 
                    \frac{1}{\w(\v)-|\w|} G(\v)F'(\v) \vec{e}_\w
                   \end{array}$    \\
      \hline
      $H^{1}_\v \times H^{1}_\v$  &   $\begin{array}{c}
                                                                  \left(\overline{F(\v)\, u_k \PB_\v}, \,\overline{G(\v)\, u_l \PB_\v}\right)\\
                                                                  \left(\overline{F(\v)\, u_k \PB_\v},\, \overline{\PB_{u_s}}\right)\\
                                                                  \left(\overline{\PB_{u_s}},\, \overline{\PB_{u_t}}\right)
                                                               \end{array}$ & 
                                                          $\begin{array}{c}
                                                                0\\
                                                                F(\v)\, u_k\PB_{u_s}\\
                                                                0
                                                             \end{array}$  \\
    \hline
  \end{tabular}
  \end{center}
\end{table}
\begin{table}[ht]
  \caption{Case $\w(\v)=|\w|$. 
                 The values of the linear map $f_2^\v$ on the elements of the bases $\bb^{i}_\v$ and $\bb^{j}_\v$ of the spaces $H_\v^i$ and $H_\v^j$, for $i,j=-1, 0, 1, 2$. 
                 In this table, $F(\v), G(\v)$ are arbitrary elements of $\F[\v]$ and $0\leq k,l\leq \mu-1$ and $1\leq s,t\leq \mu-1$.}
                 \label{tab:2}
  \begin{center}
  \renewcommand{\arraystretch}{1.5}
  \begin{tabular}{|c|c|c|}
    \hline 
    $H^i_\v \times H^j_\v$ & $(\overline{\V^i},\,\overline{\V^j})\in \bb^{i}_\v\times \bb^{j}_\v$ & $f_2^\v(\overline{\V^i},\overline{\V^j})\in \g_\v^{i+j-1}$\\
    \hline
    \hline
      $H^{-1}_\v \times H^{0}_\v$   & $\left(\overline{F(\v)}, \,\overline{G(\v)\, \vec{e}_\w}\right)$ & $0$ \\
      \hline
      $H^{-1}_\v \times H^{1}_\v$   & $\begin{array}{c}
                                                               \left(\overline{F(\v)}, \,\overline{G(\v)\, u_l \PB_\v}\right)\\
                                                                \left(\overline{F(\v)},\, \overline{\PB_{u_s}}\right)
                                                                  \end{array}$ & $\begin{array}{c}0\\F'(\v)\,u_s\end{array}$  \\
      \hline
      $H^{-1}_\v \times H^{2}_\v$    & $\left(\overline{F(\v)}, \,\overline{G(\v)\, u_l \D}\right)$ &   $0$    \\
      \hline
      $H^{0}_\v \times H^{0}_\v$     &  $\left(\overline{F(\v)\, \vec{e}_\w}, \,\overline{G(\v)\, \vec{e}_\w}\right)$ & $0$ \\
      \hline
      $H^{0}_\v \times H^{1}_\v$  &   $\begin{array}{c}
                                                               \left(\overline{F(\v)\, \vec{e}_\w}, \,\overline{G(\v)\, u_l \PB_\v}\right)\\
                                                                \left(\overline{F(\v)\, \vec{e}_\w},\, \overline{\PB_{u_s}}\right)\\
                                                                {}\\
                                                                  \end{array}$ & $\begin{array}{c}0\\ 
                                                                  \Big(\frac{\w(u_s)-|\w|}{|\w|}\frac{F(\v)-F(0)}{\v}\\
                                                                  \qquad-F'(\v)\Big)u_s\vec{e}_\w\end{array}$  \\
      \hline                                                            
      $H^{0}_\v \times H^{2}_\v$     &  $\left(\overline{F(\v)\, \vec{e}_\w}, \,\overline{G(\v)\, u_l \D}\right)$ &   $0$    \\
      \hline
      $H^{1}_\v \times H^{1}_\v$  &  $\begin{array}{c}
                                                               \left(\overline{F(\v)\, u_k \PB_\v}, \,\overline{G(\v)\, u_l \PB_\v}\right)\\
                                                                \left(\overline{F(\v)\, u_k \PB_\v},\, \overline{\PB_{u_s}}\right)\\
                                                                \left(\overline{\PB_{u_s}},\, \overline{\PB_{u_t}}\right)
                                                                  \end{array}$ & $\begin{array}{c}0\\F(\v)\,u_k\PB_{u_s}\\0 \end{array}$  \\
    \hline
  \end{tabular}
  \end{center}
\end{table}
\begin{proof}[proof of theorem \ref{prop:defosLinfini}]
One can check (by a direct computation) that the following hold:
\begin{eqn}{eq:schoutenH1}
\lb{F(\v)\, u_k \PB_\v, G(\v)\, u_l \PB_\v}_S &=& 0,\\
\lb{F(\v)\, u_k \PB_\v, \PB_{u_t}}_S &=&- \p_\v\left(F(\v)\,u_k\PB_{u_t}\right),\\
\lb{ \PB_{u_s}, \PB_{u_t}}_S &=& 0,
\end{eqn}%
for all $0\leq k,l\leq \mu-1$ and all $1\leq s,t \leq \mu-1$ and for arbitrary elements $F(\v)$ and $G(\v)$ of $\F[\v]$.
Because of (3) of proposition \ref{prp:ap}, this implies that the map $\ell_2^\v$, which is the map induced by the Schouten bracket on the cohomology $H_\v$ (and also denoted by $\LB_S$), is zero when restricted to $H_\v^1\otimes H_\v^1$.

Now,  by $\ell_1^\v=0$ and the definition (\ref{eq:f1v}) of $f_1^\v$, it is straightforward to show that the skew-symmetric graded linear map 
$f_2^\v:\bigotimes\nolimits^2 H(\g_\v,\p_\v) \to \g_\v$, defined by the tables \ref{tab:1} and \ref{tab:2}, together with $\ell_2^\v = \LB_S$, satisfy the equation $(\E_2)$. In particular, let us check this on $H^1_\v\otimes H^1_\v$. Indeed, for all $0\leq k,l\leq \mu-1$ and for arbitrary elements $F(\v)$ and $G(\v)$ of $\F[\v]$, the equation $(\E_2)$ for $\xi_1 =\overline{F(\v)\, u_k \PB_\v}$ and $\xi_2 = \overline{G(\v)\, u_l \PB_\v}$ becomes, using (\ref{eq:schoutenH1}),
 \begin{eqnarray*}
 \lefteqn{\p_\v\left(f_2^\v\left(\overline{F(\v)\, u_k \PB_\v}, \overline{G(\v)\, u_l \PB_\v}\right)\right)}\\
  &=& f_1^\v\left(\overline{\lb{F(\v)\, u_k \PB_\v, G(\v)\, u_l \PB_\v}_S} \right)-\lb{F(\v)\, u_k \PB_\v, G(\v)\, u_l \PB_\v}_S\\
  &=& 0.  
 \end{eqnarray*}
 Similarly, one also obtains $\p_\v\left(f_2^\v\left(\overline{ \PB_{u_s}}, \overline{\PB_{u_t}}\right)\right) =0$, for all $1\leq s,t \leq \mu-1$. 
 Finally, for any arbitrary element $F(\v)$ of $\F[\v]$, and for all $0\leq k\leq \mu-1$ and $1\leq t \leq \mu-1$, the identities (\ref{eq:schoutenH1}) imply that the equation $(\E_2)$ for $\xi_1 =\overline{F(\v)\, u_k \PB_\v}$ and $\xi_2 =\overline{ \PB_{u_t}}$ reads as follows
 \begin{eqnarray*}
  \lefteqn{\p_\v\left(f_2^\v\left(\overline{F(\v)\, u_k \PB_\v}, \overline{\PB_{u_t}}\right)\right)}\\
   &=& f_1^\v\left(\overline{\lb{F(\v)\, u_k \PB_\v, \PB_{u_t}}_S }\right)  -\lb{F(\v)\, u_k \PB_\v, \PB_{u_t}}_S\\
  &=&  \p_\v\left(F(\v)\,u_k\PB_{u_t}\right),
\end{eqnarray*}
 where we have used that $f_1^\v\circ \p_\v =0$. This implies that if the map $f_2^\v$ takes, on $H^1_\v\otimes H^1_\v$, the values given in the tables \ref{tab:1} and \ref{tab:2}, then the previous equations are satisfied, i.e., the equation $(\E_2)$ is satisfied on $H^1_\v\otimes H^1_\v$. From now on, we fix $f_2^\v$ to take, on $H^1_\v\otimes H^1_\v$,  the values given in the tables \ref{tab:1} and \ref{tab:2}.
 
 We have obtained the existence of the maps $\ell_1^\v$, $\ell_2^\v$ and $f_1^\v$, $f_2^\v$, satisfying the equations $(\E_1)$, $(\E_2)$ and $(\J_1)$, $(\J_2)$, $(\J_3)$.
 By the proposition \ref{prop:magic}, this implies that there exist
 skew-symmetric graded linear maps 
$$
f_3^\v : \bigotimes\nolimits^3 H(\g_\v,\p_\v) \to \g_\v \quad\hbox{ and }\quad \ell_3^\v : \bigotimes\nolimits^3 H(\g_\v,\p_\v) \to H(\g_\v,\p_\v)
$$ 
with $\deg(f_3^\v) = -2$, $\deg(\ell_3^\v) = -1$ and satisfying the equation $(\E_{3})$. Moreover, the proposition \ref{prop:magic} also says that such a map $\ell_3^\v$ necessarily satisfies the equation $(\J_{4})$.

In the equation $(\E_n)$, we denote $T_n(\Fe_n^\v, \L_{n-1}^\v; \xi_1,\dots,\xi_n)$ by $T_n^\v(\xi_1,\dots,\xi_n)$, for $n\in\N^*$ and $\xi_1,\dots,\xi_n\in H(\g,\p_\g)$, when $\Fe_n^\v$ and $\L_{n-1}^\v$ denote the elements $\Fe_n^\v:=(f_1^\v, \dots,f_n^\v)$ and $\L_{n-1}^\v:=(\ell_1^\v,\dots,\ell_{n-1}^\v)$.
By $(\E_3)$, we have $\ell_3^\v := -p\circ T^\v_3$.
Moreover, given the maps $\ell_1^\v, \ell_2^\v, f_1^\v, f_2^\v$ as previously, one can also verify that:
$$
{T^\v_3}_{\vert_{\left(H_\v^1\right)^{\otimes 3}}}=0, 
$$
so that, ${\ell_3^\v}_{\vert_{\left(H_\v^1\right)^{\otimes 3}}}={-p\circ T_3^\v}_{\vert_{\left(H_\v^1\right)^{\otimes 3}}}=0$, and the equation $(\E_3)$ is still satisfied if we choose ${f_3^\v}_{\vert_{\left(H_\v^1\right)^{\otimes 3}}}:=0$, what we do from now on. 
Let us for example show that 
$$
T^\v_3\left(\overline{F(\v)u_l\PB_\v}, \overline{\PB_{u_s}}, \overline{\PB_{u_t}}\right) = 0,
$$
for any arbitrary element $F(\v)$ of $\F[\v]$, for all $0\leq l\leq \mu-1$ and all $1\leq s, t\leq \mu-1$. First, let us point out that, by the definition (\ref{T(n)}) of $T^\v_3$, and because ${\ell^\v_2}_{\vert_{H_\v^1\otimes H_\v^1}}=0$, we simply get, for any $\xi_1,\xi_2,\xi_3\in H_\v^1$:
\begin{equation*}
T^\v_3\left(\xi_1,\xi_2,\xi_3\right) =
\lb{f^\v_1(\xi_1), f^\v_2(\xi_2,\xi_3)}_S +\lb{f^\v_2(\xi_1,\xi_2), f^\v_1(\xi_3)}_S + \lb{f^\v_2(\xi_1,\xi_3), f^\v_1(\xi_2)}_S.
\end{equation*}
Now, by the tables \ref{tab:1} and \ref{tab:2}, we obtain:
\begin{eqnarray*}
\lefteqn{T^\v_3\left(\overline{F(\v)u_l\PB_\v}, \overline{\PB_{u_s}}, \overline{\PB_{u_t}}\right) =}\\
&&\lb{F(\v)u_l\PB_{u_s}, \PB_{u_t}}_S + \lb{F(\v)u_l\PB_{u_t}, \PB_{u_s}}_S.
\end{eqnarray*}
To conclude that this is equal to zero, it suffices to show (by a direct computation) that, for any $f,g,h,l\in\A$, we have:
\begin{eqnarray*}
\lb{f\PB_l, g\PB_h}_S 
&=& f\left(\pp{l}{x}\left(\pp{g}{y}\pp{h}{z} -\pp{g}{z}\pp{h}{y}\right) + \circlearrowleft(x,y,z)\right)\D\\
&+& g\left(\pp{h}{x}\left(\pp{f}{y}\pp{l}{z} -\pp{f}{z}\pp{l}{y}\right) + \circlearrowleft(x,y,z)\right)\D\\
&=& - \lb{f\PB_h, g\PB_l}_S,
\end{eqnarray*}
where we recall that $\D$ denotes the skew-symmetric triderivation of $\A$ defined by $\D:=\pp{}{x}\wedge\pp{}{y}\wedge\pp{}{z}$, and where ``$+ \circlearrowleft(x,y,z)$" means that we consider the other terms with cyclically permuted variables $x,y,z$.

Now we have chosen the maps $\ell_1^\v$, $\ell_2^\v$, $\ell_3^\v$ and $f_1^\v$, $f_2^\v$, $f_3^\v$ such that the equations $(\E_1)$, $(\E_2)$, $(\E_3)$ and $(\J_1)$, $(\J_2)$, $(\J_3)$, $(\J_4)$ are satisfied,  ${\ell_i^\v}_{\vert_{\left(H_\v^1\right)^{\otimes i}}}=0$ for $i=2,3$, $f_2^\v$ is given by the tables \ref{tab:1} and \ref{tab:2}, and ${f_3^\v}_{\vert_{\left(H_\v^1\right)^{\otimes 3}}}=0$. 
The proposition \ref{prop:magic} once more gives us the existence of  skew-symmetric graded linear maps 
$$
f_4^\v : \bigotimes\nolimits^4 H(\g_\v,\p_\v) \to \g_\v \quad\hbox{ and }\quad \ell_4^\v : \bigotimes\nolimits^4 H(\g_\v,\p_\v) \to H(\g_\v,\p_\v)
$$ 
with $\deg(f_4^\v) = -3$ and $\deg(\ell_4^\v) = -2$ and satisfying the equation $(\E_{4})$. Moreover, according to the proposition \ref{prop:magic}, such a map $\ell_4^\v$ satisfies also the equation $(\J_{5})$. 
It is also straightforward, with the choices made previously, to show that
$$
{T^\v_4}_{\vert_{\left(H_\v^1\right)^{\otimes 4}}}=0.
$$
This implies that ${\ell_4^\v}_{\vert_{\left(H_\v^1\right)^{\otimes 4}}}= -{p\circ T^\v_4}_{\vert_{\left(H_\v^1\right)^{\otimes 4}}} =0$ and that it is possible to choose ${f_4^\v}_{\vert_{\left(H_\v^1\right)^{\otimes 4}}}=0$ (what we do from now on), so that $(\E_4)$ is still satisfied.
Finally, because ${\ell_i^\v}_{\vert_{\left(H_\v^1\right)^{\otimes i}}}=0$, for $i=2,3,4$, and ${f_i^\v}_{\vert_{\left(H_\v^1\right)^{\otimes i}}}=0$, for $i=3, 4$,  one has necessarily that:
$$
{T_{j}^\v}_{\vert_{\left(H_\v^1\right)^{\otimes j}}}=0,  \hbox{ for all } j\geq 5. 
$$
This fact, together with the proposition \ref{prop:magic}, imply that there finally exist skew-symmetric graded linear maps 
\begin{eqnarray*}
\ell_k^\v : \bigotimes\nolimits^k H(\g_\v,\p_\v) \to H(\g_\v,\p_\v), &\hbox{ with } k\geq 5,\\
f_k^\v : \bigotimes\nolimits^k H(\g_\v,\p_\v) \to \g_\v, &\hbox{ with } k\geq 5,
\end{eqnarray*}
of degrees $2-k$ and $1-k$ respectively, and satisfying 
${\ell_k^\v}_{\vert_{\left(H_\v^1\right)^{\otimes k}}}=0$, and ${f_k^\v}_{\vert_{\left(H_\v^1\right)^{\otimes k}}}=0$, for all $k\geq 5$, 
such that the maps $(\ell_1^\v,\ell_2^\v,\ell_3^\v, \dots)$ and $(f_1^\v,f_2^\v,f_3^\v,\dots)$ satisfy the conditions $(P_1)$ -- $(P_3)$, and
\begin{enumerate}
\item[-] $\left(\ell_k^\v\right)_{k\in\N^*}$ is an \l-algebra structure on $H_\v$,
\item[-] $\left(f_k^\v\right)_{k\in\N^*}$ is a quasi-isomorphism from $H_\v$ to $\g_\v$,
\end{enumerate}
hence the theorem \ref{prop:defosLinfini}.
\end{proof}
\begin{rem}
There is a natural question concerning this theorem \ref{prop:defosLinfini}, which is: is it possible that $\ell^\v_k=0$, for all $k\geq 3$? In other words, is it possible that the theorem extends to a result of \emph{formality} for $\g_\v$? Indeed, a dg Lie algebra $(\g, \p_\g, \LB_\g)$ is said to be \emph{formal} if it is linked to the dg Lie algebra $(H(\g,\p_\g), 0, \LB_\g)$ (endowed with the trivial differential and the graded Lie bracket induced by $\LB_\g$) by a quasi-isomorphism. 

In fact, we can show that, except maybe if we change the definition of $f_1^\v$ (i.e., if we consider another choice of bases $\bb^{-1}_\v, \bb^{0}_\v, \bb^{1}_\v, \bb^{2}_\v$), the map $\ell_3^\v$ cannot be zero. In the case $\w(\v)\not=|\w|$, one indeed has for example:
$$
T_3^\v\left(\bar\v,\bar\v,\bar\D\right)= 2\lb{\v, f_2^\v\left(\bar\v,\bar\D\right)}_S.
$$
We know that the choice we made for the value $f_2^\v\left(\bar\v,\bar\D\right)$ is unique, up to a $1$-cocycle for the Poisson cohomology associated to $(\A,\PB_\v)$. According to the fact that $H^1(\A,\PB_\v)\simeq\lbrace 0\rbrace$, when $\w(\v)\not=|\w|$, a $1$-cocycle is a $1$-coboundary, that is to say  an element of the form $\Vf=\Pb{F}_\v$, with $F\in\A$ (called an hamiltonian derivation). For such an element, $\lb{\v,\Vf}_S =- \Vf[\v]=0$. This implies that the value of $T_3^\v\left(\bar\v, \bar\v, \bar\D\right)$ does not depend on the choice for  $f_2^\v\left(\bar\v,\bar\D\right)$ and, using the table \ref{tab:1},
$$
T_3^\v\left(\bar\v, \bar\v, \bar\D\right)=2\lb{\v,\frac{1}{\w(\v)-|\w|}\vec{e}_\w}_S = 2\frac{\w(\v)}{|\w|-\w(\v)}\v.
$$
Because $\ell_3^\v = -p\circ T_3^\v$, we have $\ell_3^\v\left(\bar\v, \bar\v, \bar\D\right)=2\frac{\w(\v)}{|\w|-\w(\v)}\bar \v$, which is not zero. 
\end{rem}

\smallskip

\subsection{Classification of the formal deformations of $\PB_\v$}
\label{subsec:conclu}
To obtain the theorem \ref{prp:intro}, we fix an \l-algebra structure $\ell_\bullet^\v$ on $H_\v$ and a quasi-isomorphism $f_\bullet^\v$ from $H_\v$ to $\g_\v$, as in theorem \ref{prop:defosLinfini}.
By the paragraph \ref{par:formal-defos}, we know that $\Def^\nu(\g_\v)$ corresponds to the set of all the equivalence classes of the formal deformations of $\PB_\v$.
Let us now consider the set $\Def^\nu(H_\v)$. By definition of the generalized Maurer-Cartan equation (\ref{MC}) and because the \l-algebra structure $\ell_\bullet^\v=\left(\ell_k^\v\right)_{k\in\N^*}$ satisfies $\ell_1^\v=0$ and the property $(P_2)$ of the theorem \ref{prop:defosLinfini}, we have:
\begin{equation*}
\MC^\nu(H_\v) =  H_\v^1\otimes \nu\F[[\nu]] = H^2(\A,\PB_\v)\otimes \nu\F[[\nu]]. 
\end{equation*}
In the generic case, that is to say when $\w(\v) \not=|\w|$, according to proposition \ref{prp:ap}, one has $H^0_\v\simeq\lbrace0\rbrace$, so that the gauge equivalence in $\MC^\nu(H_\v)$ is trivial and
\begin{equation*}
\Def^\nu(H_\v) \simeq \MC^\nu(H_\v)=  H_\v^1\otimes \nu\F[[\nu]] = H^2(\A,\PB_\v)\otimes \nu\F[[\nu]]. 
\end{equation*}
Moreover, in the special case where $\w(\v) =|\w|$, then according to proposition \ref{prp:ap}, one has $H^0_\v = \F[\v]\overline{\vec{e}_\w}$ and in this case:
\begin{equation*}
\Def^\nu(H_\v) = H_\v^1\otimes \nu\F[[\nu]] / \sim\, = H^2(\A,\PB_\v)\otimes \nu\F[[\nu]] / \sim, 
\end{equation*}
where $\sim$ is the gauge equivalence in $\MC^\nu(H_\v)$, generated by the infinitesimal transformations of the form:
\begin{equation}\label{eq:tr}
\gamma\longmapsto  
 \gamma - \sum\limits_{k\geq 1} \frac{(-1)^{k(k-1)/2}}{(k-1) !} \ell_k^\v(\xi, \gamma,\dots,\gamma),
\end{equation}
 for $\xi = \sum\limits_{i\geq 1}\overline{ F_i(\v)\,\vec{e}_\w}\, \nu^i \in H^0_\v\otimes \nu\F[[\nu]]$, where the $F_i(\v)$ are elements of $\F[\v]$.
 We now are able to show the following:
 \begin{thm}\label{prp:conclu}
 Let $\v\in\A=\F[x,y,z]$ be a weight-homogeneous polynomial, with an isolated singularity. To $\v$ is associated the Poisson structure defined by:
 $$
 \PB_\v := \pp{\v}{x}\,\pp{}{y}\wedge\pp{}{z} + \pp{\v}{y}\,\pp{}{z}\wedge\pp{}{x} + \pp{\v}{z}\,\pp{}{x}\wedge\pp{}{y}.
 $$
 We consider the dg Lie algebra $(\g_\v,\p_\v,\LB_S)$, associated to $\v$ and defined in the paragraph \ref{subsec:jac}, of all skew-symmetric multiderivations of $\A$, equipped with the Schouten bracket $\LB_S$ and the differential $\p_\v := \lB{\PB_\v}_S$.
 
We denote by $\Cs$, the set of all $(\cc,\bar\cc)$, where
$\cc:=\left(c^k_{l,i}\in\F\right)_{{(l,i)\in\N\times\E_\v}\atop{k\in\N^*}}$ is a family of constants indexed by $\N\times\E_\v\times \N^*$
and $\bar\cc:=\left(\bar c^{\,k}_r\in\F\right)_{{1\leq r\leq \mu-1}\atop{k\in\N^*}}$ is a family of constants indexed by $\lbrace 1, \dots,Ê\mu-1\rbrace \times \N^*$, such
  that, for every $k_0\in\N^*$, the sequences $(c^{k_0}_{l,i})_{(l,i)\in\N\times\E_{\v}}$ and
  $(\bar c^{\,k_0}_r)_{1\leq r\leq\mu-1}$ have finite supports. 
  Now, for every element $(\cc, \bar\cc)=\left((c^k_{l,i}),(\bar c^{\,k}_r)\right)\in\Cs$, we associate an element $\gamma^{\cc,\bar\cc}$ of $\g_\v^1\otimes \nu\F[[\nu]]$, by the following formula:
\begin{eqnarray}\label{eq:gamma}
\gamma^{\cc,\bar\cc} := \sum_{n\in\N^*} \gamma_n^{\cc,\bar\cc} \,\nu^n,
\end{eqnarray}
with, for all $n\in\N^*$, $\gamma_n^{\cc,\bar\cc}$ given by:
\begin{eqn}[2.3]{eq:gamma-n}
 \gamma_n^{\cc,\bar\cc} &:=& \renewcommand{\arraystretch}{0.7}
 \ds\sum_{\begin{array}{c}\scriptstyle (l,i)\in\N\times\E_{\v}\\
 \scriptstyle 1\leq r\leq\mu-1\end{array}}
    \ds \sum_{\begin{array}{c}\scriptstyle a+b=n\\\scriptstyle a, b\in\N^*\end{array}} 
        c^a_{l,i}\,\bar c^{\,b}_r\, \v^l\, u_i\, \PB_{u_r}  \\
&+& \renewcommand{\arraystretch}{0.7} 
\ds\sum_{(m,j)\in\N\times\E_{\v}} c^n_{m,j}\, \v^m\,u_j\PB_\v 
\;+ \;\ds\sum_{1\leq s\leq\mu-1} \bar c^{\,n}_s\, \PB_{u_s},
\end{eqn}%
where the $u_j$, for $0\leq j\leq \mu-1$, are weight homogeneous polynomials of $\A=\F[x,y,z]$, whose images in 
$\A_{sing}(\v)=\F[x,y,z]/\ideal{\pp{\v}{x},\pp{\v}{y},\pp{\v}{z}}$ give a
basis of the $\F$-vector space $\A_{sing}(\v)$, and $u_0=1$.
 Then, one has:
 \begin{enumerate}
 \item The set of all the gauge equivalence classes of the solutions of the Maurer-Cartan equation associated to the dg Lie algebra $(\g_\v,\p_\v,\LB_S)$ is then given by:
 \begin{eqnarray*}
 \Def^\nu(\g_\v) = \lbrace \gamma^{\cc,\bar\cc}  \mid (\cc,\bar\cc)\in \Cs \rbrace /\sim,
 \end{eqnarray*}
 where $\sim$ still denotes the gauge equivalence;
 \item In the generic case where $\w(\v)\not=|\w|$, this set is exactly given by:
 \begin{eqnarray*}
 \Def^\nu(\g_\v) = \lbrace \gamma^{\cc,\bar\cc}  \mid (\cc,\bar\cc)\in \Cs \rbrace.
 \end{eqnarray*}
 \end{enumerate}
  \end{thm}

 \smallskip

 \begin{proof}
 To show this theorem, we fix an \l-algebra structure $\ell_\bullet^\v$ on $H_\v$ and a quasi-isomorphism $f_\bullet^\v$, as in theorem \ref{prop:defosLinfini}. According to the theorem \ref{thmMC}, we know that 
 $$
 \Def^\nu(\g_\v)  = \Def^\nu(f^\v_\bullet) \left(\Def^\nu(H_\v) \right).
 $$
 We also have seen at the beginning of this paragraph that, because $\ell_1^\v=0$ and because of the property $(P_2)$ of theorem \ref{prop:defosLinfini}, $\Def^\nu(H_\v) = H_\v^1\otimes \nu\F[[\nu]] / \sim$.
Now, by definition of $f_\bullet^\v$, and because it satisfies the property $(P_3)$ of the theorem \ref{prop:defosLinfini}, and by definition (\ref{eq:def(f)}) (and (\ref{eq:MCf})) of $ \Def^\nu(f^\v_\bullet)$, we have:
 \begin{eqnarray*}
 \Def^\nu(\g_\v) &=&  \Def^\nu(f^\v_\bullet) \left(\Def^\nu(H_\v) \right)\\
   &=& \left(f_1^\v + \frac{1}{2}f_2^\v \right) (H^1_\v\otimes \nu\F[[\nu]])/ \sim.
 \end{eqnarray*}
 Let $\gamma = \sum_{n\in\N^*}\gamma_n\, \nu^n$ be an element of $H^1_\v\otimes \nu\F[[\nu]]$, where each $\gamma_n$ is an element of $H^1_\v$. For $n\in\N^*$, every element $\gamma_n$ can be decomposed in the basis $\bb^1_\v$ (see the proposition \ref{prp:ap}), i.e., there exist  families of constants $(\cc, \bar\cc)=\left(\left(c^n_{m,j}\right), \left(\bar c^{\,n}_s\right)\right)\in\Cs$ satisfying:
 $$
 \gamma_n =   \renewcommand{\arraystretch}{0.7} 
\ds\sum_{(m,j)\in\N\times\E_{\v}} c^n_{m,j}\, \v^m\,u_j\PB_\v 
\;+ \;\ds\sum_{1\leq s\leq\mu-1} \bar c^{\,n}_s\, \PB_{u_s},
 $$
 for all $n\in\N$.
 Now, using the tables \ref{tab:1} and \ref{tab:2}, we obtain exactly that 
 $ \left(f_1^\v + \frac{1}{2}f_2^\v \right) \left(\gamma\right) = \gamma^{\cc,\bar\cc}$, hence the result. For the case where $\w(\v)\not=|\w|$, it only remains to recall that in this case, the gauge equivalence $\sim$ is trivial, as explained at the beginning of this paragraph.
 \end{proof}
According to what we have seen in the paragraph \ref{par:formal-defos}, the previous theorem can be translated into a result concerning the formal deformations of the family of Poisson brackets $\PB_\v$, for $\v\in\F[x,y,z]$, a weight-homogenous polynomial with an isolated singularity. It then becomes exactly the parts (a), (b) and (c) of the proposition 3.3 of \cite{AP2} and replacing $\nu\F[[\nu]]$ by $\nu\F[[\nu]]/\langle\nu^{m+1}\rangle$ (with $m\in\N^*$) in everything we have done leads to the part (d) of this proposition 3.3 of \cite{AP2}, which we write once more here:
\begin{prp}[\cite{AP2}]\label{prp:defo_A_v}
Let $\v\in\A=\F[x,y,z]$ be a weight homogeneous polynomial with an isolated
singularity. Consider the Poisson algebra $(\A,\PB_\v)$ associated
to $\v$, where $\PB_\v$ is the Poisson bracket given by (\ref{eq:PBv}). Then we have the following:
\begin{enumerate}
\item[(a)] For all families of constants 
$\left(c^k_{l,i}\in\F\right)_{{(l,i)\in\N\times\E_\v}\atop{k\in\N^*}}$ 
and $\left(\bar c^{\,k}_r\in\F\right)_{{1\leq r\leq \mu-1}\atop{k\in\N^*}}$, such
  that, for every $k_0\in\N^*$, the sequences $(c^{k_0}_{l,i})_{(l,i)\in\N\times\E_{\v}}$ and
  $(\bar c^{\,k_0}_r)_{1\leq r\leq\mu-1}$ have finite supports, the formula 
\begin{eqnarray}\label{eq:defo_qcq}
\pi_* = \PB_\v + \sum_{n\in\N^*} \pi_n \nu^n,
\end{eqnarray}
where, for all $n\in\N^*$, $\pi_n$ is given by:
\begin{eqn}[2.3]{defo_form}
 \pi_n &=& \renewcommand{\arraystretch}{0.7}
 \ds\sum_{\begin{array}{c}\scriptstyle (l,i)\in\N\times\E_{\v}\\
 \scriptstyle 1\leq r\leq\mu-1\end{array}}
    \ds \sum_{\begin{array}{c}\scriptstyle a+b=n\\\scriptstyle a, b\in\N^*\end{array}} 
        c^a_{l,i}\,\bar c^{\,b}_r\, \v^l\, u_i\, \PB_{u_r}  \\
&+& \renewcommand{\arraystretch}{0.7} 
\ds\sum_{(m,j)\in\N\times\E_{\v}} c^n_{m,j}\, \v^m\,u_j\PB_\v 
\;+ \;\ds\sum_{1\leq s\leq\mu-1} \bar c^{\,n}_s\, \PB_{u_s},
\end{eqn}%
defines a formal deformation of $\PB_\v$, where the $u_j$ ($0\leq j\leq \mu-1$) are weight homogeneous polynomials of $\A=\F[x,y,z]$, whose images in $\A_{sing}(\v)=\F[x,y,z]/\ideal{\pp{\v}{x},\pp{\v}{y},\pp{\v}{z}}$ give a basis of the $\F$-vector space $\A_{sing}(\v)$, and $u_0=1$.

\bigskip
\item[(b)] For any formal deformation $\pi_*'$ of $\PB_\v$,
there exist families of constants 
$\left(c^k_{l,i}\right)_{{(l,i)\in\N\times\E_\v}\atop{k\in\N^*}}$ 
and $\left(\bar c^{\,k}_r\right)_{{1\leq r\leq \mu-1}\atop{k\in\N^*}}$
(such that,
for every $k_0\in\N^*$, only a finite number of $c^{k_0}_{l,i}$ and
$\bar c^{\,k_0}_r$ are non-zero), 
for which $\pi_*'$ is equivalent to the
formal deformation $\pi_*$ given by the above formulas
\emph{(\ref{eq:defo_qcq})} and \emph{(\ref{defo_form})}.

\bigskip
\item[(c)] Moreover, if the (weighted) degree of the polynomial $\v$ is not equal
to the sum of the weights: $\w(\v)\not=|\w|$, then for any formal deformation $\pi_*'$ of
$\PB_\v$, there exist \emph{unique} families of constants 
$\left(c^k_{l,i}\right)_{{(l,i)\in\N\times\E_\v}\atop{k\in\N^*}}$ 
and $\left(\bar c^{\,k}_r\right)_{{1\leq r\leq \mu-1}\atop{k\in\N^*}}$
(with, for every $k_0\in\N^*$, only a finite number of non-zero $c^{k_0}_{l,i}$ and
$\bar c^{\,k_0}_r$),
such that $\pi_*'$
is equivalent to the formal deformation $\pi_*$ given by the  
formulas \emph{(\ref{eq:defo_qcq})} and \emph{(\ref{defo_form})}.

\bigskip

This means that
formulas \emph{(\ref{eq:defo_qcq})} and \emph{(\ref{defo_form})} give a \emph{system of
representatives for all formal deformations of $\PB_\v$, modulo equivalence}.

\bigskip
\item[(d)] Analogous results hold if we replace formal
  deformations by $m$-th order deformations ($m\in\N^*$) and impose in
  \emph{(c)} that $c^k_{l,i}=0$ and $\bar c^{\,k}_r=0$, as soon as $k\geq m+1$. 
\end{enumerate}
\end{prp}
\bibliographystyle{plain}
\bibliography{ref}
%


\end{document}